\documentclass[12pt,twoside,a4paper,titlepage]{amsart}
\usepackage[ansinew]{inputenc}
\usepackage{comment}
\usepackage{amsmath}
\usepackage{amsfonts}
\usepackage{amssymb}
\usepackage{mathtools}
\usepackage[arrow, matrix, curve]{xy}
\usepackage{graphicx}
\usepackage{color}
\usepackage[left=3cm,top=3cm,bottom=3.5cm]{geometry}
\usepackage{fancyhdr}
\usepackage{amsthm}
\usepackage{paralist}
\usepackage{tikz}
\usetikzlibrary{arrows}
\usetikzlibrary{decorations.markings,arrows,automata}
\usetikzlibrary{calc}
\usepackage{pgf}
\usepackage{pgfplots}
\pgfplotsset{compat=1.18}
\usepackage{textcmds}
\usepackage{microtype}
\usepackage{hyperref}
\usepackage{tabularx}
\usepackage{xy}
\usetikzlibrary{shapes.geometric}
\newcommand{\D}[1]{\textit{\textbf{#1}}}
\usepackage[nameinlink]{cleveref}
\usepackage{xcolor}
\usepackage{float}
\usepackage{multicol}
\colorlet{mygreen}{black!50!green}
\theoremstyle{definition}
\newtheorem{Def}{Definition}[section]
\newtheorem{Ex}[Def]{Example}

\theoremstyle{plain}
\newtheorem{Lem}[Def]{Lemma}
\newtheorem{Theo}[Def]{Theorem}
\newtheorem{Prop}[Def]{Proposition}
\newtheorem{Cor}[Def]{Corollary}
\theoremstyle{remark}
\newtheorem{Rem}[Def]{Remark}
\newtheoremstyle{beweis}
{3pt}
{3pt}
{default}
{\parindent}
{default}
{:}
{\newline}
{}
\theoremstyle{beweis}

\newtheoremstyle{TheoremNum}
{\topsep}{10pt}              
{\itshape}                      
{}                              
{\bfseries}                     
{.}                             
{ }                             
{\thmname{#1}\thmnote{ \bfseries #3}}
\theoremstyle{TheoremNum}
\newtheorem{thmn}{Theorem}

\newtheorem{propo}{Proposition}
\crefname{Ex}{Example}{Examples}
\crefname{Prop}{Proposition}{Propositions}
\crefname{Rem}{Remark}{Remarks}
\crefname{Cor}{Corollary}{Corollaries}
\crefname{Lem}{Lemma}{Lemmas}
\crefname{Theo}{Theorem}{Theorems}
\crefname{Def}{Definition}{Definitions}
\crefname{Not}{Notation}{Notations}

\hypersetup{
	colorlinks=true,
	linkcolor=blue,
	filecolor=magenta,      
	urlcolor=cyan,
}
\title{The space of Discrete Morse Functions}
\author{Julian Br\"uggemann}
\address[Julian Br\"uggemann]{Max Planck Institute for Mathematics, Bonn, Germany/ Mathematical Institute of the Polish Academy of Sciences, Warsaw, Poland}
\email{brueggemann@mpim-bonn.mpg.de / julian.brueggemann@impan.pl}
\usetikzlibrary{patterns}
\begin{document}
	\begin{abstract}
    In this work, we introduce a combinatorial-geometric model for the space of discrete Morse functions on any CW complex $X$. We relate this version of a space of discrete Morse functions to the space of cellular filtrations of $X$ and discuss its relationship to various concepts such as smooth Morse theory, Cerf theory, complexes of discrete Morse matchings, and induced merge trees and barcodes. 
	\end{abstract}
	\begingroup
	\let\newpage\relax
	\maketitle
	\endgroup 
	\tableofcontents
	\section{Introduction}
    Discrete Morse theory, originally introduced by Forman \cite{Forman98}, is a powerful framework for investigating the simple homotopy type and other properties of CW complexes. It has proved to be useful for the study of persistent invariants of filtered CW complexes as they appear in topological data analysis. Other than that, discrete Morse theory has numerous applications in both, pure and applied mathematics. \par
    The idea of studying combinatorial and geometric properties of the collection of mathematical objects of a certain kind originates in the concept of moduli spaces as they are used in topology and algebraic geometry. The general idea of moduli spaces is to provide a topological space $\mathcal{M}$, together with a convenient way to parametrize $\mathcal{M}$, such that points of $\mathcal{M}$ correspond to equivalence classes of certain mathematical objects one wants to investigate, and geometric features of $\mathcal{M}$ reflect relevant properties of said objects under investigations. \par
    We choose the framework of hyperplane arrangements in vector spaces of discrete functions to study the space of discrete Morse functions in a way that allows the study of various properties and equivalence classes of filtrations by combinatorial-geometric methods. \par
    We describe how Formans definition of discrete Morse functions \Cref{dMfun} gives rise to the \Cref{definitionMorsearrangement} of the Morse arrangement $\mathcal{A}(X)$, a hyperplane arrangement that subdivides the vector space of discrete functions $\mathbb{R}^X$ on a CW complex $X$ into regions, some of which represent Forman-equivalence classes of discrete Morse functions, and others do not contain discrete Morse functions at all. This way, the introduced space of discrete Morse functions allows a parametrization indexed over the face poset $F(X)$, which helps us to investigate the structure of discrete Morse functions using the combinatorial properties of the Morse arrangement $\mathcal{A}(X)$.\par 
    In this work, we investigate the space of discrete Morse functions on a given CW complex $X$ and its connection to certain other concepts in mathematics from a combinatorial-geometric viewpoint. \par 
    We provide constructions of the following spaces out of spaces of discrete Morse functions:
    \begin{enumerate} 
        \item The spaces of discrete Morse matchings on regular CW complexes.
        \item The space of merge trees.
        \item The space of barcodes.
    \end{enumerate}
    Similar work has been done in \cite{CURRY2020}, \cite{LeTi22}, and \cite{CyMi20}. Whereas \cite{CURRY2020} explores different notions of equivalences of smooth Morse functions on the sphere and the corresponding moduli spaces, both \cite{LeTi22} and \cite{CyMi20} investigate fibers of two different instances of persistence maps.\par
    The spaces of Morse functions on manifolds have been introduced in \cite{Cerf} and further investigated in \cite{HaWa} as subspaces of the spaces of smooth functions $C^\infty(M)$ on compact smooth manifolds $M$, governed by a stratification given by a certain notion of regularity for critical points and critical values. The motivation to investigate Morse functions in that way was to give a solid framework for the investigation of the pseudo-isotopy versus isotopy question and its implications for the h-cobordism theorem. Cerf already found a map that relates certain subspaces of the space of Morse functions to spaces of discrete Morse functions \Cref{CerfLemma} decades before discrete Morse functions had been invented by Forman. We extend this map to path components:\par
    \mbox{}\par
    \begin{thmn}[\ref{CorCerf}]
         It follows from \Cref{maintheosmooth} that Cerf's map from \Cref{CerfLemma} extends to a map $\eta \colon \mathcal{N}\rightarrow \mathcal{M}(M_{\mathcal{N}})$, where $M_{\mathcal{N}}$ denotes the CW decomposition of $M$ induced by any Morse function in the path component $\mathcal{N}$ of $M$. Moreover, Cerf's proof of \Cref{CerfLemma} shows that this instance of $\eta$ is a topological submersion and compatible with the respective stratifications, too.
    \end{thmn}
    
    We also give a brief description of what Cerf's combinatorial description of strata of functions looks like in the discrete case. \par
    \mbox{}\par
    Moreover, we consider the space of discrete Morse matchings, which has been introduced in \cite{ChariJo} and further investigated in \cite{CapiMini} and \cite{LIN2021102250}. We show that the space of discrete Morse matchings is also canonically associated to the Morse arrangement: \par
    \mbox{}\par
\begin{propo}[\ref{PropSpaceOfMorseMatchings}]
    Let $X$ be a CW complex, let $\mathcal{A}$ be the Morse arrangement on $\mathbb{R}^X$ and let $\mathfrak{M}(X)$ be the space of discrete Morse matchings on $X$. Then there is a canonical embedding of posets $\mathfrak{M}(X)\subset \mathcal{L}(\mathcal{A})$, where $\mathcal{L}(\mathcal{A})$ denotes the intersection poset of $\mathcal{A}$.
\end{propo}
\mbox{}\par
We describe the space of merge trees in a similar way as the space of discrete Morse functions. We model merge trees as maps from combinatorial merge trees to the real numbers and topologize their space using a mixture of the Euclidean distance and a kind of edit distance (see \Cref{definitionspaceofmergetrees}). Our edit path distance is similar to an edit distance in \cite{WeGa22}, although their edit moves focus on edges, whereas our approach focuses on nodes. \par
\mbox{}\par
\begin{propo}[\ref{PropEuclideanEditDistanceIsMetric}]
    The edit path distance $\operatorname{d}$ from \Cref{definitionspaceofmergetrees} is a pseudo-metric on the space of merge trees $Mer$ and a metric on the space of strict merge trees $Mer_<$ and the space of well-branched merge trees $Mer_{wb}$. 
\end{propo}
Moreover, we show the following:\par
\mbox{} \par
\begin{thmn}[\ref{TheoInducedMergeTreeIsContinuous}]
    Let $X$ be a regular CW complex. Then the map $M \colon \mathcal{M}(X)\rightarrow Mer$ from \Cref{inducedmergetreemap} that maps a discrete Morse function to its induced merge tree is Lipschitz continuous.
\end{thmn}
For the relationship to the space of barcodes, we use a similar approach\footnote{One can identify the space of barcodes over a given combinatorial barcode in this work with the space $\mathbb{R}^{2n}$ in \cite[Section 4.1]{BrGa}.} as in \cite{BrGa} for the space of barcodes and a construction for the barcode induced by a merge tree similar to the one in \cite{CURRY2024102031}. We obtain the following result:\par
\mbox{}\par
\begin{propo}[\ref{PropMapInducedBarcode}]
    The map $B \colon Mer_{wb} \rightarrow Bar$ from \Cref{DefInducedBarcode} is Lipschitz continuous. 
\end{propo}
This way, we provide an alternative proof for the stability of merge trees and degree 0 barcodes, in our case with respect to notions of distance that are motivated by this Morse theoretic approach. 
	\subsection*{Acknowledgements}
    The author would like to thank Max Planck Institute for Mathematics for the great scientific environment in which this project was started and the majority of the work was conducted. The author would also like to thank the Dioscuri Centre in TDA, located at the Mathematical Institute of the Polish Academy of Science, where the finishing touches were added. Moreover, the author would like to thank Andrea Bianchi, Florian Kranhold and Paul M\"ucksch for helpful discussions about the project. Most notably the author thanks his PhD advisor, Viktoriya Ozornova, for her advice, the many helpful discussions, and the detailed feedback at multiple occasions. Furthermore, the author thanks Clemens Bannwart for fruitful discussions on the topics of cell equivalences and gradient fields.

    \section{Hyperplane Arrangements and Their Associated Spaces}
    We review some well established notions from combinatorial geometry and topology, which will prove to be useful for our endeavor. The standard references we refer to are \cite{OrMa} and \cite{AguMaha}. While most notions are standard in the literature, we adapt and extend some of them to the setting of this article.
    \begin{Def}[compare \cite{OrMa} and \cite{AguMaha}]\mbox{}\\
	A \D{real hyperplane arrangement} $\mathcal{A}$ is a finite set of hyperplanes embedded in a finite-dimensional real vector space $V$. A real hyperplane arrangement is called \D{central} if all hyperplanes contain the origin. For a central arrangement $\mathcal{A}$, the \D{center} of $\mathcal{A}$ is the subspace $\mathcal{Z}(\mathcal{A})\coloneqq \bigcap_{H\in\mathcal{A}} H \subset V$ given by the intersection of all hyperplanes of $\mathcal{A}$. A central arrangement is called \D{essential} if $\mathcal{Z}(\mathcal{A}) =\{0\}$. For any central arrangement $\mathcal{A}\subset V$ we call the arrangement induced by the projection $V\rightarrow \mathcal{Z}(\mathcal{A})^{\perp}$ to the orthogonal complement of the center the \D{essential arrangement associated to} $\mathcal{A}$.\par
    An \D{oriented hyperplane arrangement} is a real hyperplane arrangement together with a choice of an orientation class of normal vectors for each hyperplane. \par
	We call the path components of the complement $V\setminus \mathcal{A}$ the \D{(open) regions} of the arrangement $\mathcal{A}$. Moreover, we call intersections of half-spaces associated to $\mathcal{A}$, with at least one half-space per hyperplane chosen, \D{faces} of the hyperplane arrangement.
    For oriented hyperplane arrangements, each of these faces $\sigma$ is associated with a sign vector $w \in \{+,-,0\}^\mathcal{A}$, where $+$ means that $\sigma$ lies on the positive side of the corresponding hyperplane, $-$ means that $\sigma$ lies on the negative side of the corresponding hyperplane, and $0$ means that $\sigma$ lies inside corresponding hyperplane. The regions correspond to sign vectors with no entry 0. The \D{face poset} $\Sigma (\mathcal{A})$ is the set of faces of $\mathcal{A}$ ordered by inclusion.
	\end{Def}
    While hyperplane arrangements can be defined in arbitrary vector spaces, it is often useful to use coordinate spaces, i.e. spaces that come with a choice of basis. This allows us to use equations in the ordered coordinates in order to specify the arrangements we consider. We consider in particular partially ordered coordinate spaces. 
    \begin{Def}
        A \D{partially ordered coordinate space} is a real vector space $V$ together with a choice of a basis that is indexed over some partially ordered set $P$. 
    \end{Def}
    For convenience, we choose a specific categorical point of view for partially ordered coordinate spaces:
    \begin{Def}\label{definitioncatofposetindvect}
    The \D{category of partially ordered coordinate spaces} $po\operatorname{Vect}$ consists of
    \begin{enumerate}
        \item[Objects:] vector spaces $\operatorname{Map}(P,\mathbb{R})$ for finite posets $P$ together with the standard basis induced by $P$, and 
        \item[Morphisms:] maps $\operatorname{Map}(Q,\mathbb{R})\rightarrow \operatorname{Map}(P,\mathbb{R})$ induced by maps $P \rightarrow Q$.
    \end{enumerate}
\end{Def}
\begin{Rem}\label{Rem:FunctorialityOfPartiallyCoordinateSpaces}
    The category of partially ordered coordinate spaces is related to the image of the contravariant functor $\operatorname{Map}(\_ , \mathbb{R})\colon \operatorname{Poset} \rightarrow \operatorname{Vect}_\mathbb{R}$, where $\operatorname{Poset}$ denotes the category of finite posets and monotone maps, and $\operatorname{Vect}_\mathbb{R}$ denotes the category of finite-dimensional real vector spaces and linear maps. In addition to the data of the image of the functor, the specific choice of partially ordered basis given by the domain poset
    is part of data of partially ordered coordinate spaces.
\end{Rem}
In this work, we consider partially ordered coordinated spaces where the bases are indexed over either face posets of CW complexes or posets associated to merge trees or barcodes. We recall the definition of the face poset of any given CW complex.
\begin{Def}
    Let $X$ be a CW complex. We say that a cell $\sigma$ is a \D{face} of a cell $\tau$ if $\sigma \subseteq \tau$. The \D{face poset} of $X$ is the set of cells of $X$, denoted by $F(X)$, together with a partial order given by $\sigma \leq\tau$ if and only if $\sigma$ is a face of $\tau$.
\end{Def}
We can characterize the three frameworks mentioned above as certain types of posets, which makes partially ordered coordinate vector spaces a unifying framework for them:
\begin{Def}\label{definitionposettypes}
    A poset $P$ is of \D{regular CW type} if it is isomorphic to the face poset $F(X)$ of a regular CW complex $X$.\par
    A poset $P$ is of \D{merge tree type} if it
    \begin{enumerate}
        \item is finite,
        \item has a unique maximum $r$, 
        \item for each $a\in P$ the set $[a,r]$ is a chain, and
        \item each non-minimal element $a\in P$ has at least two elements $b\neq b'$ such that $b\prec a$ and $b' \prec a$, where $\prec$ denotes cover relation.\footnote{The cover relation $\prec$ in posets is defined as $a\prec b$ if and only if $a< b$, and there is no $c$ such that $a< c < b$.}
    \end{enumerate}
    A poset $P$ is of \D{barcode type} if it is a finite disjoint sum of chains of length one. We call these chains the \D{bars} of the barcode. 
    \end{Def}
    \begin{Rem}
        In the definition above, we deviated from the standard convention in the literature to consider trees oriented away from the root. We do so because merge trees come with a preferred orientation from leaves to the root, induced by increasing filtration levels. \par
        In order to see that \Cref{definitionposettypes} models merge trees, we remark that (2) and (3) model arbitrary rooted trees: the chains in (3) correspond to shortest paths between the root and other nodes of the tree. The feature that inner nodes of merge trees correspond to mergers of path components, i.e.~each inner node has at least two children, is reflected in (4).  \par 
        Since merge trees can, in this fashion, be considered as both, posets and trees in the graph-theoretic sense, we use the notions from both fields interchangeably. For example, we call minimal elements leaves, non-minimal elements inner nodes, and refer to the maximum as the root.
    \end{Rem} 
    \begin{Ex}
      We consider the following example of a manifold with a height function, the induced merge tree and the corresponding poset of merge tree type.
      \begin{figure}[H]
          \centering
           \begin{tikzpicture}[scale=0.25]
     	\draw (22,1) node {};
     	\draw (0,-20) node {};
     	\begin{scope}[overlay]
     		\draw[fill=gray!50!white]  (5,0) ..  controls  (1,-1) and (1,-16) .. (1,-16) .. controls (1,-18) and (3,-21) .. (4,-16) ..  controls  (5,-8) and (6,-4) .. (8,-12) ..  controls  (9,-14) and (10,-13) .. (11,-8) ..  controls  (11,-1) and (6,0) .. (5,0) ;
     	\end{scope}
     	\draw[fill=white] (6,-4) ellipse  (2cm and 1cm);
     	\draw (6,3) node {Filtered Space:};
     	\draw (14,-10) node {$\rightarrow$};
     	\draw (14,-8.5) node {\textcolor{blue}{height}};
     	\draw[thick,->] (17,-19) -- (17,0);
     	\draw (17,1) node {$\mathbb{R}$};
     	\draw (25,3) node {Merge Tree:};
     	\draw (16+7,-19) node {$\bullet$};
     	\draw (20+7,-13) node {$\bullet$};
     	\draw (18+7,-8) node {$\bullet$};
      \draw (16+7,-20) node {$\alpha$};
     	\draw (20+7,-14) node {$\beta$};
     	\draw (18+8,-8) node {$\gamma$};
     	\draw (16+7,-19) .. controls (16.25+7,-16) and (16.5+7,-10) .. (18+7,-8);
     	\draw (20+7,-13) .. controls (19+7,-8.5) .. (18+7,-8);
     	\draw (18+7,-8) -- (18+7,0);
    \draw (40,4) node {Poset of};
    \draw (40,2) node {Merge Tree Type:};
        \node (0) at (38,-15) {$\alpha$};
        \node (1) at  (42,-15) {$\beta$};
        \node (2) at (40,-10) {$\gamma$};
        \node[rotate=60] (3) at ($(0)!0.5!(2)$) {$\leq$};
        \node[rotate=300,xscale=-1] at ($(1)!0.5!(2)$) {$\leq$};
          \end{tikzpicture}
          \caption{An example of a merge tree induced by a height function on a manifold and the corresponding poset of merge tree type.}
          \label{fig:enter-label}
      \end{figure}
    \end{Ex}
    One classic example of hyperplane arrangements, which is of particular importance for this work, is the braid arrangement. 
    \begin{Ex}\label{defbraidarrangement}
    We denote by $\mathcal{H}_n$ the hyperplane arrangement in $\mathbb{R}^n$ given by the equations $x_i=x_j$ for all $i\neq j$. We call $\mathcal{H}_n$ the \D{braid arrangement}
     $\mathcal{H}_n$ (see \cite[Example 2.3.3]{OrMa}).
    \end{Ex}

    We recall the definitions of the intersection poset and the poset of flats, which are in a sense dual to each other:
    \begin{Def}[\mbox{\cite[Definition 2.1.3]{OrMa}},\mbox{\cite[1.3.1]{AguMaha}}]
        The \D{intersection poset} $\mathcal{L}(\mathcal{A})$ of a hyperplane arrangement $\mathcal{A}$ is the set of intersections of subfamilies of $\mathcal{A}$, ordered by reverse inclusion. The \D{poset of flats} $\Pi(\mathcal{A})$ is as a set the same as the intersection poset but ordered by inclusion. The elements of the poset of flats, or the intersection poset, respectively, are called the \D{flats} of $\mathcal{A}$.
    \end{Def}
    If $\mathcal{A}$ is central, then $\mathcal{L}(\mathcal{A})$ is a geometric lattice of rank $r(\mathcal{A})$. Since all arrangements in this work are central,  we also refer to $\mathcal{L}(\mathcal{A})$ as the \D{intersection lattice} of $\mathcal{A}$. 
    It is straightforward to see that 
     $\Pi(\mathcal{A})$ is never the face poset of a simplicial complex if $\mathcal{A}$ is central and $\lvert \mathcal{A}\rvert >1$.\par
    On the other hand, the face poset $\Sigma (\mathcal{A})$ always has the structure of the face poset of a regular CW complex (see \cite[1.1.8]{AguMaha}). 
    \begin{Def}[\mbox{\cite[1.1.9]{AguMaha}}]
    \label{defsimplicialarrangement}
        We call a hyperplane arrangement simplicial if $\Sigma(\mathcal{A})\setminus \{\mathcal{Z}(\mathcal{A})\}$ has the structure of a face poset of a pure\footnote{A simplicial complex is called \D{pure} if all its inclusion-maximal simplices have the same dimension $d$.} simplicial complex. 
    \end{Def}
    In fact, $\Sigma(\mathcal{A})\setminus \{\mathcal{Z}(\mathcal{A})\}$ always gives a regular CW decomposition of the sphere in $\mathcal{Z}^\perp$. Hence, the real condition for a hyperplane arrangement being simplicial is that $\Sigma(\mathcal{A})\setminus \{\mathcal{Z}(\mathcal{A})\}$ needs to be isomorphic to the face poset of a simplicial complex. Then any such simplicial complex will be pure.
	\section{The Space of Discrete Morse Functions}
	Our goal is to present a framework which allows us to analyze geometric properties of, and homotopies between, discrete Morse functions on (regular) CW complexes. It turns out that manipulations of discrete Morse functions, like cancellations of critical cells and reordering critical cells by swapping their critical values can be geometrically realized by homotopies between discrete Morse functions, i.e.\ paths in the space of discrete Morse functions.\par
    In order to define the space of discrete Morse functions on a CW complex, we consider the space of all discrete functions, which the space of discrete Morse functions will turn out to be a subspace of. As common when working with discrete Morse theory, we adopt a slight abuse of notation by not distinguishing between a CW complex and its set of cells. 
	\begin{Def}\label{definitionspaceofdiscretefunctions}
		Let $X$ be a finite CW complex. We call the space $\mathbb{R}^X\coloneqq\mathbb{R}^{F(X)} \coloneqq \prod\limits_{\sigma \in X} \mathbb{R_\sigma}$ the \D{space of discrete functions on $X$}. 
	\end{Def}
	Here, we notice that the coordinates of the space of discrete functions are non-trivially partially ordered rather than just linearly ordered. This partial order will help us to investigate which of the discrete functions are discrete Morse functions. \par
	For that, we recall the definition of discrete Morse functions as given by Forman:
	\begin{Def}[{\cite[Definition 2.1]{Forman98}}]\label{dMfun}
		Let $X$ be a finite CW complex and let $F(X)$ be the face poset of $X$. A \D{discrete Morse function} on $X$ is a function $f\colon F(X)\rightarrow \mathbb{R}$ such that for every $p$-dimensional cell $\alpha^{(p)} \in X$ 
		\begin{enumerate}[(1)]
			\item $ \#\{ \beta^{(p+1)} \supset \alpha \vert f(\beta)\leq f(\alpha)\}\leq 1$, and 
			\item $ \#\{ \gamma^{(p-1)} \subset \alpha \vert f(\gamma)\geq f(\alpha)\}\leq 1$.
            \item $f(\alpha) <f(\beta)$ whenever $\alpha \subset \beta$ is not regular.
		\end{enumerate}
		Cells for which the inequalities in (1) and (2) both are strict are called \D{critical}. 	
	\end{Def}
    \begin{Rem}
        We include the general version for non-regular CW complexes mostly for the comparison to the smooth case in \Cref{smoothcase}. For obvious reasons, we can drop condition (3) if $X$ is regular. It is highly inconvenient to work with discrete Morse theory on non-regular CW complexes because in that case the face poset does not contain all information of the homotopy type at hand. Moreover, we make the definition for the non-regular case a bit stricter than Forman in order to simplify the definition of the induced gradient field, i.e.\ Morse functions are automatically strictly monotone along face inclusions of codimension higher than one. 
    \end{Rem}
	We recall that by \cite[Lemma 2.5]{Forman98} for each cell $\alpha \in X$ at most one of the inequalities (1) and (2) from \Cref{dMfun} may actually be an equality, whereas the other one has to be a strict inequality. Therefore, for each cell $\alpha \in X$ such that \textcolor{blue}{(1)}/\textcolor{red}{(2)} is an equality there is exactly one cell \textcolor{blue}{$\beta \supset \alpha$}/\textcolor{red}{$\gamma \subset \alpha$} such that \textcolor{blue}{$f(\beta) \leq f(\alpha)$}/\textcolor{red}{$f(\alpha)\leq f(\gamma)$}.
	Hence, any discrete Morse function $f\colon F(X) \rightarrow \mathbb{R}$ induces a partial matching $\nabla f$ on $X$ by matching $\alpha$ with \textcolor{blue}{$u(\alpha)\coloneqq \beta$}/\textcolor{red}{$d(\alpha) \coloneqq \gamma$}. The unmatched cells are called \D{critical}. To be precise, $\nabla f$ is a partial matching on the Hasse diagram $D(X)$ of the face poset $F(X)$ of $X$: 
    \begin{Def}
     Let $X$ be a CW complex. The \D{Hasse diagram} of the face poset of $X$ is the directed graph $D(X)$ that has vertices $V(D(X))\coloneqq
    F(X)$, that is the cells of $X$, and an oriented edge $(\tau,\sigma)$ whenever $\sigma \subset \tau$ is a face relation such that there exists no $\theta \in X$ with $\sigma \subset \theta \subset \tau$ ($\sigma$ is covered by $\tau$). \par
    Moreover, for a discrete Morse function $f$ on $X$ we define the \D{modified Hasse diagram} $D_f(X)$ that arises from $D(X)$ by inverting the edges matched by $f$.
    \end{Def}
	\begin{Rem}
	    In case of a regular CW complex $X$, the directed edges correspond to face relations of codimension one. For non-regular CW complexes, edges might correspond to face relations of higher codimension. This simplifies explicit constructions of discrete Morse functions in the regular case because here, all necessary considerations can be checked exclusively on the face poset. In the more general case, the face relations that correspond to all intended matched pairs would need to be checked for regularity. \par
        Furthermore, for non-regular CW complexes, the gradient paths needed for the boundary maps of the Morse chain complex are harder to spot inside the modified Hasse diagram, since they are then not just zig-zags between any matched edges and any unmatched egdges.  
	\end{Rem}
    In the literature, there is a more specialized definition of discrete Morse functions which was first proposed by Benedetti in \cite{Bene}:
       \begin{Def}[{\cite[Section 2.1]{Bene}}]\label{DefMBFunction}
	        Let $X$ be a CW complex. A \D{Morse--Benedetti} function on $X$ is a function on the face poset $f\colon F(X) \rightarrow \mathbb{R}$ that fulfills for all cells $\sigma,\tau$ in $X$ that:
	        \begin{itemize}
	       \item[\textbf{Monotonicity: }] If $\sigma \subset \tau$ we have $f(\sigma)\leq f(\tau)$.
           \item[\textbf{Semi-injectivity: }] $\lvert f^{-1}(\{z\})\rvert \leq 2$ for all $z\in \mathbb{R}$.
           \item[\textbf{Genericity: }] If $f(\sigma)=f(\tau)$, then either $\sigma \subseteq \tau$ or $\tau \subseteq \sigma$ holds.
           \item[\textbf{Regularity: }] If $\sigma \subset \tau$ is irregular, then $f(\sigma) < f(\tau)$.
	        \end{itemize}
	\end{Def}
    It is straightforward to see that Morse--Benedetti functions are always discrete Morse functions but in general not the other way around.\par
    For some purposes, it is useful to weaken Benedetti's definition:
    \begin{Def}\label{DefWeakMBFunction}
        Let $X$ be a CW complex. We call a function $f\colon F(X) \rightarrow \mathbb{R}$ a \D{weak Morse--Benedetti function} if it is a discrete Morse function and satisfies semi-injectivity and genericity in the sense of \Cref{DefMBFunction}.
    \end{Def}
    The only difference between weak Morse--Benedetti functions and Morse--Benedetti functions is that matched cells of Morse--Benedetti functions have the same value, whereas for weak Morse--Benedetti functions, the higher-dimensional cell might have a strictly smaller value than the smaller-dimensional cell of the matched pair.\par
    \begin{Ex}
        We consider the following example of three discrete Morse function that induce the same matching. 
        \begin{figure}[H]
            \centering
            \begin{tikzpicture}
                \draw (0,10) -- node[midway,above] {1} (2,10) -- node[midway,right] {2} (2,8) --  node[midway,below] {5} (0,8) -- node[midway,left] {5} (0,10);
                \draw[color=red,->] (2,10) -- (1,10);
                \draw[color=red,->] (2,8) -- (2,9);
                \draw (0,10) node {$\bullet$};
                \draw[color=red] (2,10) node {$\bullet$};
                \draw (0,8) node {$\bullet$};
                \draw[color=red] (2,8) node {$\bullet$};
                
                \draw (-0.25,10.25) node {0};
                \draw (2.25,10.25) node {1};
                \draw (-0.25,7.75) node {3};
                \draw (2.25,7.75) node {3};
            \end{tikzpicture}
            \hspace{1cm}
            \begin{tikzpicture}
                \draw (0,10) -- node[midway,above] {1} (2,10) -- node[midway,right] {2} (2,8) --  node[midway,below] {5} (0,8) -- node[midway,left] {6} (0,10);
                \draw[color=red,->] (2,10) -- (1,10);
                \draw[color=red,->] (2,8) -- (2,9);
                \draw (0,10) node {$\bullet$};
                \draw[color=red] (2,10) node {$\bullet$};
                \draw (0,8) node {$\bullet$};
                \draw[color=red] (2,8) node {$\bullet$};
                
                \draw (-0.25,10.25) node {0};
                \draw (2.25,10.25) node {1};
                \draw (-0.25,7.75) node {4};
                \draw (2.25,7.75) node {3};
            \end{tikzpicture}
            \hspace{1cm}
            \begin{tikzpicture}
                \draw (0,10) -- node[midway,above] {1} (2,10) -- node[midway,right] {2} (2,8) --  node[midway,below] {5} (0,8) -- node[midway,left] {4} (0,10);
                \draw[color=red,->] (2,10) -- (1,10);
                \draw[color=red,->] (2,8) -- (2,9);
                \draw (0,10) node {$\bullet$};
                \draw[color=red] (2,10) node {$\bullet$};
                \draw (0,8) node {$\bullet$};
                \draw[color=red] (2,8) node {$\bullet$};
                
                \draw (-0.25,10.25) node {0};
                \draw (2.25,10.25) node {1};
                \draw (-0.25,7.75) node {3};
                \draw (2.25,7.75) node {2};
            \end{tikzpicture}
            \caption{Three different dMfs that induce the same matching (indicated in red). The one in the middle is weakly Morse--Benedetti but not Morse--Benedetti. The one on the right is Morse--Benedetti. }
            \label{fig:Ex:DMfvswMBvsMB}
        \end{figure}
        The dMf on the left is not weakly Morse--Benedetti due to the two critical edges labeled 5. The dMf in the middle is not Morse--Benedetti because of the matched pair with labels 2 and 3.
    \end{Ex}
	The definition of discrete Morse functions gives rise to a subdivision of the space of discrete functions induced by an oriented hyperplane arrangement in $\mathbb{R}^X$ as follows: \par 
	 
	
	We observe that each point $p\in\mathbb{R}^X$ represents a discrete function $f\colon F(X)\rightarrow \mathbb{R}$ by $f(\sigma)=p_\sigma$. The conditions (1) and (2) (and (3)) of \Cref{dMfun} impose conditions on the components of a discrete function $p\in \mathbb{R}^X$ whenever there is a face relation of codimension 1 (or a cover relation) between the corresponding cells. That is, for the question whether $p$ represents a discrete Morse function it is relevant whether the inequality $p_\sigma < p_\tau$ holds whenever $\sigma \subset \tau$ is a face relation of codimension 1 (or a cover relation). This information is equivalent to the question on which side of the hyperplane $H_\sigma^\tau$ in $\mathbb{R}^X$ given by $x_\sigma=x_\tau$ the point $p$ lies. Hence, the conditions (1) and (2) (and (3)) of \Cref{dMfun} can be checked by analyzing the position of $p$ relative to the hyperplane arrangement $\mathcal{A}$ given by the equations $x_\sigma=x_\tau$ for all cover face relations\footnote{Face relations of codimension 1 are a special case of cover relations.} $\sigma \subset \tau$ in $X$. We define the following oriented hyperplane arrangement:
	\begin{Def}\label{definitionMorsearrangement}
		Let $E$ be the set of cover face relations in $X$. We define the hyperplane arrangement $\mathcal{A}(X)\coloneqq \{H_\sigma^\tau \vert \sigma\subset \tau \in E\}$ where $H_\sigma^\tau$ is given by the equation $x_\sigma=x_\tau$. We observe that the assignment $(\sigma\subset \tau)\mapsto H_\sigma^\tau$ defines a bijection $E\cong \mathcal{A}$. Then let $\mathcal{S}(X)\coloneqq\{+,-,0\}^E$  be the  sign vectors for the oriented hyperplane arrangement $\mathcal{A}^+(X)$ obtained by orienting $\mathcal{A}(X)$ in the following way: For each hyperplane $H_\sigma^\tau$ we define the open half-space defined by $x_\sigma < x_\tau$ to be the positive side of  $H_\sigma^\tau$, the open half-space given by $x_\sigma > x_\tau$ to be the negative side of $H_\sigma^\tau$, and on $H_\sigma^\tau$ itself we have the value 0. 
  We call $\mathcal{A}^+(X)$ the \D{Morse arrangement} on $\mathbb{R}^X$, respectively on $X$.\par
		For any sign vector $x\in \mathcal{S}(X)$ we define the sets $x_+\coloneqq \{e\in E \vert x_e=+\}$, $x_-\coloneqq \{e\in E \vert x_e=-\}$, and $x_0\coloneqq \{e\in E \vert x_e=0\}$.
	\end{Def}
\begin{Rem}
	The index set $E$ of $\mathcal{S}(X)$ is by construction the set of edges of $D(X)$. Hence, we will use the two different points of view interchangeably without explicitly mentioning it in the notation from now on. Moreover, we drop the $+$ in the notation since there is only one orientation we are interested in. Furthermore, we will denote $\mathcal{A}(X)$ and $\mathcal{S}(X)$ just by $\mathcal{A},\mathcal{S}$, respectively, when there are not multiple complexes to be confused with.
\end{Rem}
	We interpret an element $x\in \mathcal{S}$ as an equivalence class of discrete functions on a given regular CW complex $X$.
    For two functions $p,q\in \mathbb{R}^X$ we define $p \sim q$ if and only if for all hyperplanes $H\in \mathcal{A}(X)$, the function $p$ lies on the positive side of $H$ if and only if $q$ lies on the positive side of $H$. It turns out that this way the regions that contain discrete Morse functions can be identified with their Forman equivalence classes, that is, the ones that induce the same acyclic matching. 
    The positive region, i.e.\ the region on the positive side of all hyperplanes, turns out to be exactly the set of all critical discrete Morse functions on $X$. A point $p\in \mathbb{R}^X$ being inside a hyperplane $H_\sigma^\tau$ or on the negative side of  $H_\sigma^\tau$ indicates that if $p$ represents a discrete Morse function, then the cells $\sigma$ and $\tau$ are matched by $p$. Moreover, a function $p\in \mathbb{R}^X$ is a discrete Morse function if and only if all functions in the same region as $p$ are. \par
    It is immediate that the Morse arrangement $\mathcal{A}\subset \mathcal{H}_{\lvert X \rvert}$ is a subarrangement of the braid arrangement \Cref{defbraidarrangement}. Moreover, $\mathcal{A}$ is always central but never essential, i.e.\ the intersection of all hyperplanes of $\mathcal{A}$ always contains $0\in \mathbb{R}^X$ but never consists of just the origin. \par
    We give a characterization of discrete Morse functions in the language of hyperplane arrangements:
    \begin{Def}\label{DefinitionMorseFunctionHyperplane}
        Let $X$ be a CW complex and let $\mathcal{A}$ be the corresponding Morse arrangement. A discrete function $p\in \mathbb{R}^X$ is a discrete Morse function if for every cell $\sigma \in X$
        \begin{enumerate}[(1)]
			\item $ \#\{ \tau \in X \vert (\tau,\sigma)\in D(X),\ p_\tau\leq p_\sigma \}\leq 1$, and 
			\item $ \#\{ \tau \in X \vert (\sigma,\tau)\in D(X),\ p_\sigma\leq p_\tau \}\leq 1$.
            \item $p_\sigma < p_\tau$ whenever $(\sigma,\tau) \in D(X)$ corresponds to a non-regular cover relation. 
		\end{enumerate}
        We call any region of $\mathcal{A}$ that consists of discrete Morse functions a \D{Morse region}. We call the union of all Morse regions of $\mathcal{A}$ the \D{space of discrete Morse functions} on $X$ and denote it by $\mathcal{M}(X)$. We call the unique region $R_{cr}(X)$ that only contains critical discrete Morse functions, i.e.~discrete Morse functions for which every cell is critical, the \D{critical region} of $\mathcal{A}$. Moreover, we define the \D{space of discrete Morse filtrations}\footnote{Lying in the closure of the critical region ensures the discrete function to be weakly increasing. This, in turn, ensures that for every cell the function value coincides with its filtration value in the sublevel filtration.} $\mathcal{M}_{fil}(X)$ as the intersection $\mathcal{M}_{fil}(X)\coloneqq \overline{R_{cr}(X)}\cap \mathcal{M}(X)$. We define the \D{space of Morse--Benedetti functions} $\mathcal{MB}(X)\subset \mathcal{M}(X)$ as the subspace of those discrete Morse functions that satisfy monotonicity, semi-injectivity and genericity. We define the \D{space of weak Morse--Benedetti funtcions} $\mathcal{MB}_w(X)\subset \mathcal{M}(X)$ as the subspace of those discrete Morse functions that satisfy semi-injectivity and genericity.
    \end{Def}
    \begin{Rem}
    The Morse arrangement is in general not a reflection arrangement and not simplicial. For example, the critical region of the Morse arrangement on the 2-simplex $\Delta^2$ has 9 faces but it would need to have only 6 faces in order to be simplicial due to the dimension of $\mathcal{Z}(\mathcal{A}(\Delta^2))^\perp$. Hence, the regular CW complex associated to the Morse arrangement, as well as its subcomplex induced by Morse regions, is in general not simplicial.\par
    The reflection arrangement associated to $\mathcal{A}(X)$\footnote{I.e.~ the smallest reflection arrangement that contains $\mathcal{A}(X)$.} for a connected regular CW complex $X$ is the braid arrangement associated to the symmetric group $\mathfrak{S}_{\lvert X \rvert}$. For arbitrary CW complexes, the story is less well behaved because one has a priori no combinatorial control over cover relations being regular.
\end{Rem}
    \begin{Ex}
    We consider the space of discrete Morse functions on the interval $I$:
       \begin{figure}[H]
       \centering
       \begin{tikzpicture}[scale=1.2]
          \node[inner sep=2pt, circle] (0) at (4.5,6) [draw] {};
          \node[inner sep=2pt, circle] (1) at (5.5,6) [draw] {};
          \draw (0) -- (1);
          \draw (4.5,5.7) node {$\alpha$};
          \draw (5,5.7) node {$\gamma$};
          \draw (5.5,5.7) node {$\beta$};
          \end{tikzpicture}\\
          \begin{tikzpicture}
          \begin{axis}[
          axis x line =center,
          axis y line = center,
       xmin=-5,xmax=5,
       ymin=-5,ymax=5,
       xlabel={\tiny$\alpha-\beta$},
       ylabel={\tiny$-\alpha+2\gamma-\beta$},
       ticks=none,
       grid=both,
       ]
       \addplot[mark=none,domain=-5:5, color=blue] {x/3};
       \addplot[mark=none,domain=-5:5, color=green] {-x/3};
       \node at (axis cs:4,1.8) {\textcolor{blue}{$\alpha=\gamma$}};
       \node at (axis cs:-4,1.8) {\textcolor{green}{$\beta=\gamma$}};
       
       \node at (axis cs:1,3) {\textcolor{blue}{$+$}};
       \node at (axis cs:-1,3) {\textcolor{green}{$+$}};
        
       \node at (axis cs:1,-3) {\textcolor{blue}{$-$}};
       \node at (axis cs:-1,-3) {\textcolor{green}{$-$}};
      
       \node at (axis cs:-3,0.5) {\textcolor{blue}{$+$}};
       \node at (axis cs:-3,-0.5) {\textcolor{green}{$-$}};
       
       \node at (axis cs:3,0.5) {\textcolor{blue}{$-$}};
       \node at (axis cs:3,-0.5) {\textcolor{green}{$+$}};
       
       \node at (axis cs: 0,2) {$\bullet$};
       \node at (axis cs: 0.7,2.4) {$\operatorname{dim}$};
    \end{axis}

        \end{tikzpicture}
                   \caption{The Morse arrangement in the space of discrete functions on $I$}
    \end{figure}
    In order to draw the picture we made use of the fact that summands of multiples of the vector $(1,1,\dots ,1)$ have no effect on whether a point in $\mathbb{R}^X$ corresponds to a Morse region or not. Hence, we only drew the orthogonal complement $\operatorname{span}\langle(1,1,\dots ,1)\rangle^\perp$. That is, we actually drew the essential arrangement associated to $\mathcal{A}$.
    \end{Ex}
    \begin{Prop}
    Let $X$ be a regular CW complex, then the space of discrete Morse functions $\mathcal{M}(X)$ is contractible.
    \end{Prop}
    \begin{proof}
        We prove the slightly stronger statement that $\mathcal{M}(X)$ is a star domain with the dimension function $p_\sigma=\operatorname{dim}(\sigma)$ as a star point. Let $f$ be a discrete Morse function, we show that $f_t\coloneqq (1-t)\cdot f +t\cdot \operatorname{dim}$ is a Morse function for all $t\in [0,1]$. In order to do that, we prove that for any face relation $\sigma \subset \tau$ of codimension one where $f_t(\tau) \leq f_t(\sigma)$ holds, the inequality $f(\tau) \leq f(\sigma)$ must also hold. For $t=1$ we have $f_1=\operatorname{dim}$ and the statement is true. For $t<1$ we have:
        \begin{align*}
            f(\tau) = \frac{f_t(\tau)-t\cdot(\operatorname{dim}(\sigma)+1)}{1-t} \leq \frac{f_t(\tau)-t\cdot \operatorname{dim}(\sigma)}{1-t}\leq \frac{f_t(\sigma)-t\cdot \operatorname{dim} (\sigma)}{1-t}=f(\sigma).
        \end{align*}
    This implies that $f_t$ is a discrete Morse function for all $t$ if $f$ is a discrete Morse function. Hence, we connected any discrete Morse function to the dimension function by a line segment of discrete Morse functions. In particular, $\mathcal{M}(X)$ is contractible.
    
    \end{proof}
    \begin{Prop}
        Let $X$ be a regular CW complex. Then we have $\mathcal{MB}(X)\subset \overline{R_{cr}(X)}$.
        Moreover, $\mathcal{MB}_w(X)$ is dense in $\mathcal{M}(X)$ and $\mathcal{MB}(X)$ is dense in $\mathcal{M}_{fil}(X)$.
    \end{Prop}
    \begin{proof}
        Since Morse--Benedetti functions\footnote{See \Cref{DefMBFunction}} are monotone, they can only lie on the positive side of or on a hyperplane $H\in \mathcal{A}(X)$ but never on the negative side. Hence, all Morse--Benedetti functions lie in the closure of the critical region $\overline{R_{cr}(X)}$. \par
        Due to genericity, weak Morse--Benedetti functions cannot have the same value on different cells unless such cells are faces of one another. By semi-injectivity, this can only happen if the two cells are faces of codimension one, i.e.~if they are matched. Hence, discrete Morse functions can only fail to be weak Morse--Benedetti functions on intersections of at least one hyperplane of $\mathcal{H}\setminus \mathcal{A}$. Thus, $\mathcal{M}(X)\setminus \mathcal{MB}_w(X)$ is contained in the first stratum of the stratification induced by $\mathcal{H}(X)$ which means that $\mathcal{MB}_w(X)$ is dense in $\mathcal{M}(X)$. It follows that $\mathcal{MB}(X)$ is dense in $\mathcal{M}_{fil}(X)$ by restricting the stratification\footnote{Any hyperplane arrangement induces a stratification of the ambient vector space where the $k$th-stratum is given by the union of all intersections of $k$ hyperplanes.} induced by $\mathcal{H}$ to $\mathcal{M}_{fil}(X)$.
    \end{proof}
    Having established that the space of discrete Morse functions is homotopically trivial, we are going to investigate its combinatorial and geometrical properties. Since properties of objects are ideally investigated via morphisms, we take a look at meaningful notions of morphisms in this context. The goal is to identify morphisms between CW complexes that induce a meaningful notion of morphisms between spaces of discrete Morse functions. Such a notion of morphisms $\varphi \colon X\rightarrow Y$ should
    \begin{enumerate}
        \item induce maps on the spaces of discrete functions,
        \item preserve regularity of face relations, and
        \item induce maps of hyperplane arrangements that are compatible with the orientation.         
    \end{enumerate}
    We address these points in the order mentioned above, starting with a notion of morphism that induces maps on spaces of discrete functions.
    \begin{Def}\label{DefCellMapOrderPreservingMapAndSimplicialMap}
        Let $X,Y$ be CW complexes. \par
        A continuous map $\varphi\colon X \rightarrow Y$ is called \D{cellular} if it maps the $n$-skeleton $X_n$ to the $n$-skeleton $Y_n$. \par
        A continuous map $\varphi\colon X \rightarrow Y$ is called \D{order-preserving} if it induces an order-preserving map $F(\varphi)\colon F(X)\rightarrow F(Y)$, i.e.~if for any cell $\sigma \in X$, there is a cell $\sigma' \in Y$, such that $f(\sigma) \subseteq \sigma'$ and for cells $\sigma\subset\tau\in X $ we have $f(\sigma) \subset f(\tau)\in Y$. \par
        We call a map $\varphi\colon X \rightarrow Y$ \D{simplicial} if it is cellular and order preserving.
    \end{Def}
    \begin{Rem}\label{RemAbstractVSGeometric}
        One central point of this work is the interplay between CW complexes and their face posets. Although discrete Morse theory excels on simplicial complexes, we need to include CW complexes for the comparison to smooth Morse theory in \Cref{smoothcase}. Hence, on one hand we need to treat discrete Morse functions on CW complexes and on the other hand we want to treat discrete Morse functions on abstract simplicial complexes. We consider both points of view on regular CW complexes as a middle ground, since homotopy types of regular CW complexes are uniquely determined by their face posets. Nonetheless, it is important to point out that even though the face poset determines its corresponding regular CW complex up to homotopy, geometric cellular morphisms of regular CW complexes are more general than morphisms of the corresponding face posets. \par 
        We also remark that in \Cref{DefCellMapOrderPreservingMapAndSimplicialMap} cellular maps belong to the framework of CW complexes, whereas order-preserving maps are closely related to face posets. We introduce a notion of simplicial maps as a middle ground between the two.
    \end{Rem}
    \begin{Ex}
    In order to see that cellular maps do not necessarily induce maps on face posets, consider a regular CW decomposition of the circle $S^1$ in two 0-cells $a,b$ and two 1-cells $A,B$. Let a continuous map $f\colon S^1\rightarrow S^1$ be given by $f(a)=f(b)=f(A)=a$ and $f$ is defined on $B$ by winding around the circle twice. Then $f$ does not induce a map on face posets because the image $f(B)$ contains multiple cells, in fact all cells of $S^1$. Nonetheless, $f$ maps the 0-skeleton to the 0-skeleton and the 1-skeleton to the 1-skeleton. Since $f$ winds around the circle twice, the map $f$ is also not homotopic to any map that could induce a map on face posets.  
    \end{Ex}
    Having established notions of morphisms that induce maps on face posets, and thus induce maps of spaces of discrete functions in the sense of partially ordered coordinate spaces, we continue with notions of morphisms that preserve regularity of face relations. In order to do so, we focus on two approaches, which we then combine for regular CW complexes. The first approach uses the combinatorial structure of the faces in the boundary of cells. The second one uses the notion of cell equivalences, which combines homotopy and cell structure. 
    \begin{Def}\label{Def:Non-DegenerateMap}
        Let $X$ and $Y$ be CW complexes. We say that a cell $\tau \in X$ \D{is spanned by} a set of 0-cells $\{\sigma_i\}\subset X$ if $\{\sigma_i\}= (\partial \tau)_0$. A simplicial map $\varphi\colon X \rightarrow Y$ is called \D{non-degenerate} if for any set of 0-cells $\{\sigma_i\}\subset X$ such that  an $n$-cell $\sigma$ is spanned by $\{\sigma_i\}\subset X$ for some $n$, the minimal $k$-cell $\tau$ such that $\varphi(\sigma)\subseteq \tau$ for some $k\leq n$ is spanned by the set $\{\varphi(\sigma_i)\}\subset Y$.
    \end{Def}
\begin{Ex}\label{ExDegAndNonDegMaps}
We consider some examples and non-examples for non-degenerate maps:
\begin{enumerate}
    \item Inclusion of subcomplexes $X\hookrightarrow Y$ are non-degenerate.
    \item Let $A\subset X$ be a subcomplex. Then the quotient map $X \rightarrow X/A$ that collapses $A$ to a point is non-degenerate.
    \item Let $P\rightarrow Q$ be a poset map. Then the induced map on order complexes\footnote{Also known as the classifying space of the poset.} $\Delta(X)\rightarrow\Delta(Y)$ is non-degenerate. 
    \item Let $I$ be the standard interval considered as a cell complex with two 0-cells that are connected by a 1-cell. 
    $$\begin{tikzpicture}
        \draw (0,0) node {$\bullet$};
        \draw (2,0) node {$\bullet$};
        \draw (0,0) -- (2,0);
        \draw (0,-0.3) node {$0$};
        \draw (2,-0.3) node {$1$};
        \draw (1,-0.3) node {$\frac{1}{2}$};
    \end{tikzpicture}$$
    Consider the map 
    $$f\colon I \rightarrow I, f(t)\coloneqq \begin{cases}
        t & t\leq \frac{1}{2}\\
        1-t & t> \frac{1}{2}
    \end{cases}. $$
    Then $f$ is simplicial and even preserves the cover relation in $F(I)$, but $f$ is not non-degenerate.\par
    Similar degenerate collapse maps of a cell to one of its faces also exist in higher dimensions. 
\end{enumerate}
\end{Ex}
In the literature, there is an established notion of equivalence of CW complexes that takes both, combinatorial and topological features, into account: 
    \begin{Def}[\mbox{\cite[Definition 2.1]{FRANKS1979199},~\cite[Definition 6.32]{banyaga2004lectures}}]\label{cell equivalence}
        Let $X$ and $X'$ be finite CW complexes. For any cell $e\subset X$ denote by $X(e)$ the smallest subcomplex of $X$ containing $e$. \par
        The complexes $X$ and $X'$ are called \D{cell equivalent} if and only if there is a homotopy equivalence $h\colon X \rightarrow X'$ with the property that
        $h$ induces an isomorphism of face posets\footnote{In \cite{FRANKS1979199}, the author originally mentioned a ``one-to-on correspondence" of cells. We interpret the term ``correspondence" in the more precise above-mentioned way of an induced map on face posets, which also seems to be what the author of \cite{FRANKS1979199} intended.} $F(h)\colon F(X) \xrightarrow{\cong} F(X')$
        such that if $e \subset X$ corresponds to $e'\subset X'$, then $h$ maps $X(e)$ to $X'(e')$ and restricts to a homotopy equivalence of these subcomplexes.
    \end{Def}
        
    
    \begin{Rem}\label{RemarkaboutCellEquivalence}
    It follows directly from the definition that cell equivalences are in particular simplicial and homotopy equivalences. For the context of discrete Morse theory, we prefer the point of view of simplicial maps between regular CW complexes because they focus more on the combinatorial structure of the face poset. However, in the context of arbitrary CW complexes one has to work with cell equivalences because they keep track of the irregularity of attaching maps. \par
    Furthermore, cell equivalences induce isomorphisms of the corresponding face posets. This way, cell equivalences and non-degenerate maps are the link between geometric and abstract morphisms of regular CW complexes.
    \end{Rem}    

    \begin{Lem}\label{nondegmapsaresimplicialandinjective}
        Let $\varphi \colon X \rightarrow Y$ be a simplicial map of regular CW complexes.
        If $\varphi$ preserves the cover relation $\prec$ in the corresponding face posets, then $\varphi$ maps $n$-cells to $n$-cells.
    \end{Lem}
    \begin{proof}
        We prove the statement by induction: Since simplicial maps are cellular, they map 0-cells to 0-cells. Let $\sigma$ be an $n$ cell of $X$. By inductive assumption, all $n-1$ cells $\sigma' \subset \partial \sigma$ are mapped to $n-1$ cells of $Y$. Because $\sigma'\prec \sigma$ and $\varphi$ preserves cover relations, we have $\varphi (\sigma') \prec \varphi (\sigma)$. Because $Y$ is regular, this implies $\operatorname{dim} \varphi(\sigma)=n$. 
    \end{proof}
    \begin{Prop}\label{PropNonDegenerateCellEquiv}
        Let $\varphi \colon X \rightarrow Y$ be a non-degenerate map of regular CW complexes that preserves the cover relation. Then $\varphi$ induces homotopy equivalences on all subcomplexes $X(\sigma)$ for $\sigma$ any cell of $X$.\par
        Moreover, if the induced map $F(\varphi)$ is injective, then $\varphi$ is a cell equivalence onto its image.
    \end{Prop}
    \begin{proof}
        We observe that for regular CW complexes $X$, the subcomplexes $X(\sigma)$ for an $n$-cell $\sigma$ are homeomorphic to $D^n$. 
        Since $\varphi$ preserves cover relations, it also maps chains of cover relations to chains of cover relations. Since $\varphi$ is non-degenerative, for every cell $\sigma$ it preserves the initial vertices of all maximal chains of cover relations that terminate in $\sigma$. We prove inductively that for each such chain of cover relations $\sigma_0 \prec \dots \prec \sigma_i \prec \dots \prec \sigma$ in $X$ that finishes in $\sigma \in X$, the restriction of $\varphi$ to each $\sigma_i$ is a cell equivalence onto its image.\par
        For $i=0$, there is nothing to show. For $i=1$, it follows from preservation of cover relations and non-degeneracy, that any 1-cell $\sigma_1$ that is spanned by two 0-cells $\alpha,\beta$ needs to be mapped to a 1-cell that is spanned by $\varphi(\alpha)$ and $\varphi(\beta)$. 
        Moreover, $\varphi$ is injective on the 0-cells of $\sigma$, because otherwise either 
        \begin{enumerate}
            \item a 1-cell between two 0-cells would be mapped to a 0 cell, which would violate preservation of cover relations, or 
            \item such a 1-cell would be mapped to a degenerate 1-cell, i.e. both 0-cells of $\sigma_1$ would be mapped to the same 0-cell. 
        \end{enumerate}
        Thus, $\varphi_{\lvert X(\sigma)}$ is a continuous extension of a bijection $S^0 \cong S^0$, and, hence, must be a homotopy equivalence. \par
        Let $\sigma_i$ be an $i$-cell in the boundary of $\sigma$. By inductive assumption, $\varphi_{\lvert \partial \sigma_i}$ is a cell equivalence. In particular, $\varphi (\partial \sigma_i)\simeq S^{n-1}$. Thus, preservation of cover relations and non-degeneracy imply that $\sigma_i$ is mapped to a cell $\tau \in Y$ such that the gluing map of $\tau$ is a homeomorphism onto $\varphi (\partial \sigma_i)$, i.e.~the map $\varphi_{\lvert X(\sigma_i)}$ is a continuous extension of a cell equivalence $\partial \sigma_i \simeq \partial \tau$. This implies that $\varphi_{\lvert X(\sigma)}$ is a homotopy equivalence because $Y(\tau)\cong D^n$ is contractible. \par
        For the second statement, we observe that if $F(\varphi)$ is injective, then
        we know that the induced map $F(\varphi) \colon F(X)\rightarrow F(\operatorname{im}(\varphi))$ is an isomorphism of face posets.
        Together with the first statement, it follows directly from \Cref{cell equivalence} that $\varphi \colon X \rightarrow \varphi (X)$ is a cell equivalence.
    \end{proof}
    \begin{Rem}
    It turns out that all assumptions of \Cref{PropNonDegenerateCellEquiv} are necessary: If $\varphi$ is assumed to be non-degenerate but does not necessarily preserve cover relations, then $\varphi$ might be just a quotient map to a point. \par
    If $\varphi$ is assumed to preserve cover relations but not necessarily to be non-degenerate, it might be a map as in \Cref{ExDegAndNonDegMaps} (4).\par 
    In both of these cases, $\varphi$ would not induce homotopy equivalences on subcomplexes of the form $X(\sigma)$.\par
    Moreover, if we assume $\varphi$ to induce homotopies on subcomplexes of the form $X(\sigma$ but do not necessarily assume injectivity of the induced map on face posets, $\varphi$ could identify pairs of parallel cells as in the following example:
    $$
    \begin{tikzpicture}
        \draw[green] (2,0) arc (0:180:1);
        \draw[green] (2,0) arc (360:180:1);
        \draw[color=blue] (0,0) node {$\bullet$};
        \draw[color=blue] (2,0) node {$\bullet$};
        \draw (3,0) node {$\rightarrow$};
        \draw[color=green] (4,0) -- (6,0);
        \draw[color=blue] (4,0) node {$\bullet$};
        \draw[color=blue] (6,0) node {$\bullet$};
    \end{tikzpicture}
    $$
    In particular, the map $\varphi$ would not need to be a cell equivalence.
    \end{Rem}
    Having established the relationship between non-degenerate maps and cell equivalences, we want to use these notions to induce maps on spaces of discrete Morse functions. The straightforward way to do that is to use the contravariant functoriality of partially ordered coordinate spaces, see \Cref{Rem:FunctorialityOfPartiallyCoordinateSpaces}.
    \begin{Prop}\label{proposition: pullback map}
        Let $\varphi \colon X \rightarrow Y$ be a non-degenerate map of CW complexes. Then there is an induced linear map $\varphi^*\colon \mathbb{R}^Y \rightarrow \mathbb{R}^X$ given by $\varphi^*g\coloneqq g\circ \varphi$, for $g\colon Y \rightarrow \mathbb{R}$.\par
        If $\varphi$ is injective, then $\varphi^{*{-1}}(\mathcal{A}(X))\subset \mathcal{A}(Y)$ is a subarrangement. Moreover, $\varphi^*$ restricts in the injective case to a map $\mathcal{M}(Y)\rightarrow\mathcal{M}(X)$.
    \end{Prop}
    \begin{proof}
        The first statement is true because non-degenerate maps induce maps on face posets.
        The second statement follows because non-degenerate maps preserve cover relations. \par
        For the third statement we consider that injectivity implies that each pair of cells of $X$ that are related by a cover relation gets their function values from distinct cover relations in $Y$. Hence, the numbers from (1) and (2) of \Cref{DefinitionMorseFunctionHyperplane} are preserved by $\varphi^*$.
    \end{proof}
    \begin{Def}
        We call the induced map $\varphi^*$ from \Cref{proposition: pullback map} the \D{pullback map} induced by $\varphi$.
    \end{Def}
    Since we require the relatively strong notion of non-degenerate maps, we also get a map in the other direction.
    \begin{Lem}\label{embeddingofMorsearrangement}
        Let $\varphi\colon X \rightarrow Y$ be a non-degenerate map between CW complexes. Then $\varphi$ induces a linear map $\varphi_*\colon \mathbb{R}^X \rightarrow \mathbb{R}^Y$ given by
        $$\varphi_*(f)(\sigma) \coloneqq \begin{cases}
            f(\varphi^{-1}(\sigma)) & \text{ if }\lvert \varphi^{-1}(\sigma) \rvert =1 \\
            0 & else
        \end{cases} $$
        that maps $\mathcal{A}^+(X)$ to $\mathcal{A}^+(Y)$. If $\varphi$ is injective, then $\varphi_*$ embeds $\mathcal{A}^+(X)$ into $\mathcal{A}^+(Y)$.
    \end{Lem}
    \begin{proof}
    The first statement follows directly from the fact that $\varphi$ induces a map on face posets that preserves cover relations. The second statement follows because injectivity of $\varphi$ implies an injective correspondence of cover relations.

    \end{proof}
    \begin{Rem}
        The term non-degenerate comes from the fact that the same construction with degenerate simplicial maps might lead to mapping discrete Morse functions to discrete functions with degenerate cells.
     \end{Rem}
     In fact, if $\varphi$ is a degenerate simplicial map, $F(\varphi)$ does not preserve cover relations. Since $F(\varphi)$ is a poset map between finite posets, the fiber $F(\varphi)^{-1}(y)$ of a point $y \in F(Y)$ is either empty or an interval in $F(X)$. Intervals in $F(X)$ correspond to subspaces in $\mathbb{R}^X$ which means that $\varphi_*$ collapses the subspace corresponding to $F(\varphi)^{-1}(y)$ to one coordinate corresponding to $y$. Geometrically, this corresponds to collapsing\footnote{in the sense of collapsing a topological space, not necessarily in the sense of simple collapses} the corresponding subcomplex $U \subset X$ to a point.\par
     In order to define an induced map on the space of discrete functions one needs to make a choice for a linear map $\tilde{\varphi}_* \colon \mathbb{R}^U \rightarrow \mathbb{R}^{y}$.\par
     Canonical options like projections to one of the coordinates or taking the sum of the coordinates in general do not preserve the property of being a discrete Morse function. The only canonical possibilities that induce a map on the space of discrete Morse functions are the ones which ensure that the images of the collapsed subcomplexes become critical. If one wants to preserve matched pair of simplices along the collapsed subcomplex, one needs to apply more elaborate constructions that depend on the given input Morse function.
     \begin{Def}\label{definitionfaceonecofaceone}
         Let $X$ be a CW complex. For $\sigma \in X$, we define 
         \begin{itemize}
             \item $\operatorname{Face}(\sigma)\coloneqq \{\tau \in X\rvert \ \tau \subset \sigma \text{ arbitrary face relation}\}$,
             \item $\operatorname{Coface}(\sigma)\coloneqq \{\tau \in X \rvert \ \tau \supset \sigma \text{ arbitrary face relation}\}$,
             \item $\operatorname{Face}_1(\sigma)\coloneqq \{\tau \in X \rvert \ \tau \subset \sigma \text{ cover relation in }F(X)\}$, and
             \item $\operatorname{Coface}_1(\sigma)\coloneqq \{\tau \in X \rvert \ \tau \supset \sigma \text{ cover relation in }F(X)\}$.
         \end{itemize}
     \end{Def}
     \begin{Def}\label{inducedmaponspaceofdmfs}
         Let $\varphi\colon X \rightarrow Y$ be
         a non-degenerate map between CW complexes. We define an induced map $\varphi_*\colon \mathbb{R}^X \rightarrow \mathbb{R}^Y$ on the spaces of discrete functions. We define $\varphi_*(f)$ for any discrete function $f\in \mathbb{R}^X$ as follows: 
         \begin{enumerate}
             \item for every $\tau \in Y$ such that $\lvert \varphi^{-1}(\tau)\rvert =1$, i.e.\ $\varphi^{-1}(\tau) =\{ \sigma \}$, we define $$\varphi(f)(\tau)\coloneqq f(\sigma)$$ for any $f\in \mathbb{R}^X.$
             \item for any $\tau \in Y$ such that $\lvert \varphi^{-1}(\tau)\rvert \neq 1$ we proceed as follows: \par 
             We define 
             $$
                 \operatorname{up}(f_{\lvert \operatorname{Coface}_1(\tau)}) \coloneqq
                 \begin{cases}
                \operatorname{min}(f_{\lvert \operatorname{Coface}_1(\tau)}) & \text{ if } f \text{ is defined on some cell of } \operatorname{Coface}_1(\tau). 
                 \\
                 \operatorname{min}(f_{\lvert \operatorname{Coface}(\tau)}) & \text{ if } f \text{ is only defined on cells of } 
                 \operatorname{Coface}(\tau) \setminus \operatorname{Coface}_1(\tau)\\
                  \operatorname{max}(f_{\lvert \operatorname{Face}_1(\tau)}) +2 & \text{ if } f \text{ is not defined on any cell of }\operatorname{Coface}(\tau).
                 \end{cases}
                 $$
                 Here, we slightly abuse notation because $f_{\lvert \operatorname{Coface}_1(\tau)}$ might not be defined for all cells of $\operatorname{Coface}_1(\tau)$. If $f_{\lvert \operatorname{Coface}_1(\tau)}$ is only defined for some cells of $\operatorname{Coface}_1(\tau)$, we take the minimum on those cells.
             Start with the minimal dimensional $\tau$ such that $\lvert \varphi^{-1}(\tau)\rvert \neq 1$.
             \begin{itemize}
                 \item If $\operatorname{Face}_1(\tau)=\emptyset$ and $\operatorname{Coface}_1(\tau)= \emptyset$, then define $$\varphi_*(f)_\tau\coloneqq 0.$$
                 \item If $\operatorname{Face}_1(\tau)=\emptyset$ and $\operatorname{Coface}_1(\tau)\neq \emptyset$, then define 
                 $$\varphi_*(f)_\tau \coloneqq \begin{cases}
                     
                 \operatorname{up}(f_{\lvert \operatorname{Coface}_1(\tau)}) -1 & \text{ if } f \text{is defined on any cell of } \operatorname{Coface}(\tau)\\
                 0 & \text{else}
                 \end{cases}.$$
                 
                 \item If $\operatorname{Face}_1(\tau)\neq\emptyset$ and $\operatorname{Coface}_1(\tau)= \emptyset$, then define 
                 $$\varphi_*(f)_\tau \coloneqq \operatorname{max}(f_{\lvert \operatorname{Face}_1(\tau)}) +1.$$
                 \item If $\operatorname{Face}_1(\tau)\neq\emptyset$ and $\operatorname{Coface}_1(\tau)\neq \emptyset$, then define $$\varphi_*(f)_\tau \coloneqq \frac{\operatorname{max}(f_{\lvert \operatorname{Face}_1(\tau)})+\operatorname{up}(f_{\lvert \operatorname{Coface}_1(\tau)})}{2}.$$
                 
             \end{itemize}
             We apply this definition inductively over the dimension of the cells for which $\varphi_*$ is not yet defined.
         \end{enumerate}
         We call $\varphi_*$ the \D{push forward map} induced by $\varphi$.
     \end{Def}
     \begin{Rem}
         The set $\operatorname{Face}_1(\sigma)$ corresponds to the set of elements $\tau$ of $F(X)$ which are covered by $\sigma$, whereas the cells of $\operatorname{Coface}_1(\sigma)$ are exactly the cells $\tau$ of $X$ which cover $\sigma$. If one uses the covering relation instead of the codimension, there is a straightforward generalization to non-regular CW complexes.
     \end{Rem}
     \begin{Prop}
         The induced map $\varphi_*$ is affine-linear. Moreover, $\varphi_*$ is linear if and only if $\operatorname{Face}_1(\tau)=\emptyset$ and $\operatorname{Coface}_1(\tau)= \emptyset$ for all $\tau$ such that $\lvert \varphi^{-1}(\tau)\rvert \neq 1$. In particular, $\varphi_*$ is always linear if $\varphi$ is injective.
     \end{Prop}
     \begin{proof}
         The proof is straightforward from the definition.
     \end{proof}
     \begin{Prop}
         The map $\varphi_*$ restricts to a map of spaces of discrete Morse functions $\varphi_*\colon\mathcal{M(X)} \rightarrow \mathcal{M}(Y)$.
         Moreover, $\varphi_*$ preserves induced matchings if $\varphi$ is injective.
     \end{Prop}
     \begin{proof}
         If $\varphi$ in injective, then $X$ can be interpreted as a subcomplex of $Y$. The map $\varphi_*(f)$ is then just an extension of $f$ to $Y$ that is critical on $Y\setminus \varphi(X)$. By construction, we have $\varphi_*(f)_{\lvert \varphi(X)}=f\circ \varphi^{-1}$ and, therefore, the second statement holds.\par
         For cells $\tau$ such that $\lvert \varphi^{-1}(\tau)\rvert \neq 1$, we observe that according to \Cref{inducedmaponspaceofdmfs} $\varphi_*(f)$ is constructed such that every face that is covered by $\tau$ gets a strictly smaller function value than $\tau$ and every coface that covers $\tau$ gets a strictly greater function value than $\tau$ under $\varphi_*(f)$. Hence, every $\tau$ such that $\lvert \varphi^{-1}(\tau)\rvert \neq 1$ is by construction a critical cell of $\varphi_*(f)$. In particular, $\varphi_*(f)$ is a discrete Morse function.
     \end{proof}
    \begin{Ex}\label{Exampleinduceddmf}
    We consider the following non-degenerate map between two regular complexes:
    \begin{figure}[H]
        \centering
        \begin{tikzpicture}
            \draw[fill=cyan] (5.5,2) -- (4.5,2) -- (3,0) -- (7,0) --  cycle;
            \draw[fill=lightgray] (5.5,2) -- (8,1.5) -- (7,0) -- (5.5,2) -- cycle;
            \draw[fill=lightgray] (4.5,2) -- (3,0) -- (2,1.5) -- (4.5,2) -- cycle;
            \node[inner sep=-2pt,color=green] (0) at (4.5,2) {$\bullet$};
            \node[inner sep=-2pt,color=green] (00) at (5.5,2) {$\bullet$};
            \node[inner sep=-2pt, color=purple] (1) at (3,0) {$\bullet$};
            \node[inner sep=-2pt, color=purple] (2) at (7,0) {$\bullet$};
            \node[inner sep=-2pt] (3) at (8,1.5) {$\bullet$};
            \node[inner sep=-2pt] (4) at (2,1.5) {$\bullet$};
            \draw[color=green] (0) -- (00);
            \draw[color=blue] (0) -- (1);
            \draw[color=blue] (00) -- (2);
            \draw (00) -- (3);
            \draw (2) -- (3);
            \draw (0) -- (4);
            \draw (1) -- (4);
            \draw[color=purple] (1) -- (2);
            \node[inner sep=-2pt,color=green] (0) at (4.5,2) {$\bullet$};
            \node[inner sep=-2pt,color=green] (00) at (5.5,2) {$\bullet$};
            \node[inner sep=-2pt, color=purple] (1) at (3,0) {$\bullet$};
            \node[inner sep=-2pt, color=purple] (2) at (7,0) {$\bullet$};
            \node[inner sep=-2pt] (3) at (8,1.5) {$\bullet$};
            \node[inner sep=-2pt] (4) at (2,1.5) {$\bullet$};
            \node[inner sep=-2pt] (5) at (9,1) {$\rightarrow$};
            \draw (9,1) node[above] {$\varphi$};

            \draw[fill=lightgray] (5+7,2) -- (7+7,1) -- (5+7,0) -- (5+7,2) -- cycle;
            \draw[fill=lightgray] (5+7,2) -- (5+7,0) -- (2+8,1) -- (5+7,2) -- cycle;
            \node[inner sep=-2pt,color=green] (0) at (5+7,2) {$\bullet$};
            \node[inner sep=-2pt, color=purple] (1) at (5+7,0) {$\bullet$};
            \node[inner sep=-2pt, color=purple] (2) at (5+7,0) {$\bullet$};
            \node[inner sep=-2pt] (3) at (7+7,1) {$\bullet$};
            \node[inner sep=-2pt] (4) at (2+8,1) {$\bullet$};
            \draw[color=blue] (0) -- (1);
            \draw[color=blue] (0) -- (2);
            \draw (0) -- (3);
            \draw (2) -- (3);
            \draw (0) -- (4);
            \draw (1) -- (4);
            \node[inner sep=-2pt,color=green] (0) at (5+7,2) {$\bullet$};
            \node[inner sep=-2pt, color=purple] (1) at (5+7,0) {$\bullet$};
            \node[inner sep=-2pt, color=purple] (2) at (5+7,0) {$\bullet$};
            \node[inner sep=-2pt] (3) at (7+7,1) {$\bullet$};
            \node[inner sep=-2pt] (4) at (2+8,1) {$\bullet$};
            
        \end{tikzpicture}
        \caption{A simplicial map between regular CW complexes}
    \end{figure}
    We consider the following discrete Morse function $f$:
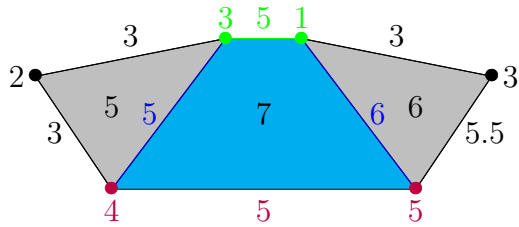
    \begin{figure}[H]
        \centering
        \begin{tikzpicture}
            \draw[fill=cyan] (5.5,2) -- (4.5,2) -- (3,0) -- (7,0) --  cycle;
            \draw[fill=lightgray] (5.5,2) -- (8,1.5) -- (7,0) -- (5.5,2) -- cycle;
            \draw[fill=lightgray] (4.5,2) -- (3,0) -- (2,1.5) -- (4.5,2) -- cycle;
            \node[inner sep=-2pt,color=green] (0) at (4.5,2) {$\bullet$};
            \node[inner sep=-2pt,color=green] (00) at (5.5,2) {$\bullet$};
            \node[inner sep=-2pt, color=purple] (1) at (3,0) {$\bullet$};
            \node[inner sep=-2pt, color=purple] (2) at (7,0) {$\bullet$};
            \node[inner sep=-2pt] (3) at (8,1.5) {$\bullet$};
            \node[inner sep=-2pt] (4) at (2,1.5) {$\bullet$};
            \draw[color=green] (0) -- (00) node[above,midway] {5};
            \draw[color=blue] (0) -- (1) node[left,midway] {5};
            \draw[color=blue] (00) -- (2) node[right,midway] {6};
            \draw (00) -- (3) node[above,midway] {3};
            \draw (2) -- (3) node[right,midway] {5.5};
            \draw (0) -- (4) node[above,midway] {3};
            \draw (1) -- (4) node[left,midway] {3};
            \draw[color=purple] (1) -- (2) node[below,midway] {5};
            \node[inner sep=-2pt,color=green] (0) at (4.5,2) {$\bullet$};
            \node[inner sep=-2pt,color=green] (00) at (5.5,2) {$\bullet$};
            \node[inner sep=-2pt, color=purple] (1) at (3,0) {$\bullet$};
            \node[inner sep=-2pt, color=purple] (2) at (7,0) {$\bullet$};
            \node[inner sep=-2pt] (3) at (8,1.5) {$\bullet$};
            \node[inner sep=-2pt] (4) at (2,1.5) {$\bullet$};
            \draw[color=green] (4.5,2) node[above] {3};
            \draw[color=green] (5.5,2) node[above] {1};
            \draw (8,1.5) node[right] {3};
            \draw (2,1.5) node[left] {2};
            \draw[color=purple] (3,0) node[below] {4};
            \draw[color=purple] (7,0) node[below] {5};
            \draw (3,1.1) node {5};
            \draw (5,1) node {7};
            \draw (7,1.1) node {6};
        \end{tikzpicture}
    \end{figure}
    The induced discrete Morse function $\varphi_*(f)$ is:
    \begin{figure}[H]
        \centering
        \begin{tikzpicture}
            \draw[fill=lightgray] (5+7,2) -- (7+7,1) -- (5+7,0) -- (5+7,2) -- cycle;
            \draw[fill=lightgray] (5+7,2) -- (5+7,0) -- (2+8,1) -- (5+7,2) -- cycle;
            \node[inner sep=-2pt, color=green] (0) at (5+7,2) {$\bullet$};
            \node[inner sep=-2pt, color=purple] (1) at (5+7,0) {$\bullet$};
            \node[inner sep=-2pt, color=purple] (2) at (5+7,0) {$\bullet$};
            \node[inner sep=-2pt] (3) at (7+7,1) {$\bullet$};
            \node[inner sep=-2pt] (4) at (2+8,1) {$\bullet$};
            \draw[color=blue] (0) -- (1) node[midway,left] {3.5};
            \draw[color=blue] (0) -- (2);
            \draw (0) -- (3) node[midway,above] {3};
            \draw (2) -- (3) node[midway,below] {5.5};
            \draw (0) -- (4) node[midway,above] {3};
            \draw (1) -- (4) node[midway,below] {3};
            \node[inner sep=-2pt] (0) at (5+7,2) {$\bullet$};
            \node[inner sep=-2pt, color=purple] (1) at (5+7,0) {$\bullet$};
            \node[inner sep=-2pt, color=purple] (2) at (5+7,0) {$\bullet$};
            \node[inner sep=-2pt] (3) at (7+7,1) {$\bullet$};
            \node[inner sep=-2pt] (4) at (2+8,1) {$\bullet$};
            \draw[color=green] (5+7,2) node [above] {2};
            \draw (5+7,0) node [color=purple,below] {2};
            \draw (7+7,1) node [right] {3};
            \draw (2+8,1) node [left] {2};
            \draw (11,1) node {5};
            \draw (13,1) node {6};
        \end{tikzpicture}
    \end{figure}
    \end{Ex}
    \begin{Rem}
        As we see in \Cref{Exampleinduceddmf}, the map $\varphi_*$ only preserves matchings where it is injective. Therefore, $\varphi_*(f)$ might in general induce a Morse complex with more cells than $f$ does even if $\varphi$ collapses a subcomplex of $X$.
    \end{Rem}
   As a particular important class of maps of CW complexes that induce maps on the spaces of discrete Morse functions, we consider internal collapses and extensions. The case of simple extensions is already covered by injective non-degenerate maps because simple extensions induce inclusions. The case of simple collapses of free faces is a bit more difficult because a pair of a free face and its unique coface in general does not induce a canonical simplicial map. In this case, we only have a (up to homotopy) canonical cellular homotopy equivalence. This still allows us to construct restricted versions of the push forward and the pullback map:
   \begin{Def}
       Let $\varphi\colon X\rightarrow X\setminus \{\sigma,\tau\}$ be a simple collapse, where $\sigma$ is the free face of $\tau$ that is being collapsed. Then we define the \D{push forward map} induced by $\varphi$ by $\varphi_* \colon \mathcal{M}(X)_{\vert x_\sigma\geq x_\tau} \rightarrow \mathcal{M}(X\setminus \{\sigma,\tau\}), \  \varphi_*(f)(\xi)\coloneqq f(\xi)$ for $\xi \in X\setminus \{\sigma,\tau\}$ and where $\mathcal{M}(X)_{\vert x_\sigma\geq x_\tau}$ is the component of $\mathcal{M}(X)$ where $\sigma$ and $\tau$ are matched.\par
       We define the \D{pullback map} induced by $\varphi$ by $\varphi^* \colon \mathcal{M}(X\setminus \{\sigma,\tau\})\rightarrow \mathcal{M}(X)_{\vert x_\sigma\geq x_\tau}$, $\varphi^*(f)(\xi)\coloneqq f(\xi)$ for $\xi \in X\setminus \{\sigma,\tau\}$, $\varphi^*(f)(\sigma)=\varphi^*(f)(\tau)\coloneqq \operatorname{max}(f_{\lvert \operatorname{Face}_1(\tau)}) +1$.
   \end{Def}
   \begin{Ex}
   We consider the following example of a simple collapse.
   \begin{center}
       \begin{tikzpicture}[decoration={markings,mark=at position 0.6 with {\arrow{triangle 60}}}]
           \draw[fill=lightgray] (0,0) -- (2,-2) -- (2,2) -- cycle;
           
           \node[circle, fill=black, inner sep=1.5pt, outer sep=0pt] (0) at (0,0) {};
           \node[circle, fill=black, inner sep=1.5pt, outer sep=0pt] (1) at (2,-2) {};
           \node[circle, fill=black, inner sep=1.5pt, outer sep=0pt] (2) at (2,2) {};

           \node[anchor = west] at (2,0) {\small{$\sigma$}};
           \node at (1,0) {\small{$\tau$}};
           \draw[->,color=red] (2,0) -- (1.3,0);

            \node at (3,0) {$\rightarrow$};
            \node at (3,0.3) {$\varphi$};
            \node[circle, fill=black, inner sep=1.5pt, outer sep=0pt] (a) at (0+4,0) {};
           \node[circle, fill=black, inner sep=1.5pt, outer sep=0pt] (b) at (2+4,-2) {};
           \node[circle, fill=black, inner sep=1.5pt, outer sep=0pt] (c) at (2+4,2) {};
           \draw (a) -- (b);
           \draw (a) -- (c);
       \end{tikzpicture}
       \end{center}
       We illustrate the pullback and push forward map in the two separate examples below.
       \begin{center}\begin{minipage}{0.48\textwidth}
           \begin{tikzpicture}
               \draw[fill=lightgray] (0,0) -- (2,-2) -- (2,2) -- cycle;
           
           \node[circle, fill=black, inner sep=1.5pt, outer sep=0pt] (0) at (0,0) {};
           \node[circle, fill=black, inner sep=1.5pt, outer sep=0pt] (1) at (2,-2) {};
           \node[circle, fill=black, inner sep=1.5pt, outer sep=0pt] (2) at (2,2) {};

            \node[anchor = east] at (0) {0};
            \node[anchor = south east] at (1,1) {1};
            \node[anchor = south] at (2) {1};
           \node[anchor = west] at (2,0) {\small{$4$}};
           \node at (1,0) {\small{$3$}};
           \node[anchor = north east] at (1,-1) {1};
           \node[anchor = north] at (1) {0};

            \node at (3,0) {$\rightarrow$};
            \node at (3,0.3) {$\varphi_*$};
            \node[circle, fill=black, inner sep=1.5pt, outer sep=0pt] (a) at (0+4,0) {};
           \node[circle, fill=black, inner sep=1.5pt, outer sep=0pt] (b) at (2+4,-2) {};
           \node[circle, fill=black, inner sep=1.5pt, outer sep=0pt] (c) at (2+4,2) {};
           \draw (a) -- (b);
           \draw (a) -- (c);

           \node[anchor = east] at (a) {0};
            \node[anchor = south east] at (1+4,1) {1};
            \node[anchor = south] at (c) {1};
           \node[anchor = north east] at (1+4,-1) {1};
           \node[anchor = north] at (b) {0};
           \end{tikzpicture}
       \end{minipage}
           \begin{minipage}{0.48\textwidth}
           \begin{tikzpicture}
               \draw[fill=lightgray] (0,0) -- (2,-2) -- (2,2) -- cycle;
           
           \node[circle, fill=black, inner sep=1.5pt, outer sep=0pt] (0) at (0,0) {};
           \node[circle, fill=black, inner sep=1.5pt, outer sep=0pt] (1) at (2,-2) {};
           \node[circle, fill=black, inner sep=1.5pt, outer sep=0pt] (2) at (2,2) {};

\node[anchor = east] at (0) {0};
            \node[anchor = south east] at (1,1) {4};
            \node[anchor = south] at (2) {3};
           \node[anchor = west] at (2,0) {\small{$5$}};
           \node at (1,0) {\small{$5$}};
           \node[anchor = north east] at (1,-1) {2};
           \node[anchor = north] at (1) {1};

            \node at (3,0) {$\leftarrow$};
            \node at (3,0.3) {$\varphi^*$};
            \node[circle, fill=black, inner sep=1.5pt, outer sep=0pt] (a) at (0+4,0) {};
           \node[circle, fill=black, inner sep=1.5pt, outer sep=0pt] (b) at (2+4,-2) {};
           \node[circle, fill=black, inner sep=1.5pt, outer sep=0pt] (c) at (2+4,2) {};
           \draw (a) -- (b);
           \draw (a) -- (c);

           \node[anchor = east] at (a) {0};
            \node[anchor = south east] at (1+4,1) {4};
            \node[anchor = south] at (c) {3};
           \node[anchor = north east] at (1+4,-1) {2};
           \node[anchor = north] at (b) {1};
           \end{tikzpicture}
       \end{minipage}
       \end{center}
   \end{Ex}
   We extend these constructions to internal collapses and extensions.
   \begin{Def}\label{definition: pullback and push forward for internal collapses and extensions}
       Let $\varphi\colon X\rightarrow Y $ be an internal collapse that collapses a pair of cells $(\sigma,\tau)$ and let $\psi\colon Y\rightarrow X$ be the inverse internal extension.\par
       We define the \D{push forward map} $\varphi_*\colon \mathcal{M}(X)_{\vert x_\sigma\geq x_\tau} \rightarrow \mathcal{M}(Y)$ by $\varphi_*(f)(\xi)\coloneqq f(\xi) \text{ for }\xi \neq \sigma,\tau$\par
       We define the \D{pullback map} as follows:
       $$\varphi^*(f)(\xi)\coloneqq \begin{cases} f(\xi) & \text{for }\xi \neq \sigma,\tau \\
       f(\sigma)=f(\tau)\coloneqq \frac{\operatorname{max}(f_{\lvert \operatorname{Face}_1(\tau)})+\operatorname{up}(f_{\lvert \operatorname{Coface}_1(\sigma)})}{2}
       \end{cases},
       $$
       where $\operatorname{up}$ is defined as in \Cref{inducedmaponspaceofdmfs} and $\operatorname{Coface}_1$ and $\operatorname{Face}_1$ are defined as in \Cref{definitionfaceonecofaceone}.\par
       For the internal extensions, we define the push forward/pullback as the pullback/push forward of the inverse internal collapse.
   \end{Def}
   \begin{Ex}We consider the following example of an internal collapse.
       \begin{center}
       \begin{tikzpicture}[decoration={markings,mark=at position 0.6 with {\arrow{triangle 60}}}]
           \draw[fill=lightgray] (0,0) -- (2,-2) -- (4,0) -- (2,2) -- cycle;
           
           \node[circle, fill=black, inner sep=1.5pt, outer sep=0pt] (0) at (0,0) {};
           \node[circle, fill=black, inner sep=1.5pt, outer sep=0pt] (1) at (2,-2) {};
           \node[circle, fill=black, inner sep=1.5pt, outer sep=0pt] (2) at (2,2) {};
            \node[circle, fill=black, inner sep=1.5pt, outer sep=0pt] (3) at (4,0) {};
           \draw (1) -- (2);
           \node[anchor = west] at (2,0) {\small{$\sigma$}};
           \node at (1,0) {\small{$\tau$}};
           \draw[->,color=red] (2,0) -- (1.3,0);

            \node at (5,0) {$\rightarrow$};
            \node at (5,0.3) {$\varphi$};
            \draw[fill=lightgray] (0+6,0) -- (2+6,-2) -- (4+6,0) -- (2+6,2) -- cycle;
            \node[circle, fill=black, inner sep=1.5pt, outer sep=0pt] (a) at (0+6,0) {};
           \node[circle, fill=black, inner sep=1.5pt, outer sep=0pt] (b) at (2+6,-2) {};
           \node[circle, fill=black, inner sep=1.5pt, outer sep=0pt] (c) at (2+6,2) {};
           \node[circle, fill=black, inner sep=1.5pt, outer sep=0pt] (d) at (4+6,0) {};
           \draw (a) -- (b);
           \draw (a) -- (c);
           \draw (d) -- (b);
           \draw (d) -- (c);
       \end{tikzpicture}
       \end{center}
       One instance for the push forward map could be:
        \begin{center}
       \begin{tikzpicture}[decoration={markings,mark=at position 0.6 with {\arrow{triangle 60}}}]
           \draw[fill=lightgray] (0,0) -- (2,-2) -- (4,0) -- (2,2) -- cycle;
           
           \node[circle, fill=black, inner sep=1.5pt, outer sep=0pt] (0) at (0,0) {};
           \node[circle, fill=black, inner sep=1.5pt, outer sep=0pt] (1) at (2,-2) {};
           \node[circle, fill=black, inner sep=1.5pt, outer sep=0pt] (2) at (2,2) {};
            \node[circle, fill=black, inner sep=1.5pt, outer sep=0pt] (3) at (4,0) {};
           \draw (1) -- (2);
           \node[anchor = east] at (0) {0};
           \node[anchor = south east] at (1,1) {4};
            \node[anchor = south] at (2) {3};
           \node[anchor = west] at (2,0) {\small{$7$}};
           \node at (1,0) {\small{$5$}};
           \node[anchor = north east] at (1,-1) {2};
           \node[anchor = north] at (1) {1};
           \node[anchor = west] at (3) {5};
           \node[anchor = south west] at (3,1) {6};
           \node[anchor = north west] at (3,-1) {6};
           \node at (3,0) {8};

            \node at (5,0) {$\rightarrow$};
            \node at (5,0.3) {$\varphi_*$};
            \draw[fill=lightgray] (0+6,0) -- (2+6,-2) -- (4+6,0) -- (2+6,2) -- cycle;
            \node[circle, fill=black, inner sep=1.5pt, outer sep=0pt] (a) at (0+6,0) {};
           \node[circle, fill=black, inner sep=1.5pt, outer sep=0pt] (b) at (2+6,-2) {};
           \node[circle, fill=black, inner sep=1.5pt, outer sep=0pt] (c) at (2+6,2) {};
           \node[circle, fill=black, inner sep=1.5pt, outer sep=0pt] (d) at (4+6,0) {};
           \draw (a) -- (b);
           \draw (a) -- (c);
           \draw (d) -- (b);
           \draw (d) -- (c);

           \node[anchor = east] at (a) {0};
           \node[anchor = south east] at (1+6,1) {4};
            \node[anchor = south] at (c) {3};
           \node[anchor = north east] at (1+6,-1) {2};
           \node[anchor = north] at (b) {1};
           \node[anchor = west] at (d) {5};
           \node[anchor = south west] at (3+6,1) {6};
           \node[anchor = north west] at (3+6,-1) {6};
           \node at (2+6,0) {8};
       \end{tikzpicture}
       \end{center}
       And one instance for the pullback map could be:
       \begin{center}
       \begin{tikzpicture}[decoration={markings,mark=at position 0.6 with {\arrow{triangle 60}}}]
           \draw[fill=lightgray] (0,0) -- (2,-2) -- (4,0) -- (2,2) -- cycle;
           
           \node[circle, fill=black, inner sep=1.5pt, outer sep=0pt] (0) at (0,0) {};
           \node[circle, fill=black, inner sep=1.5pt, outer sep=0pt] (1) at (2,-2) {};
           \node[circle, fill=black, inner sep=1.5pt, outer sep=0pt] (2) at (2,2) {};
            \node[circle, fill=black, inner sep=1.5pt, outer sep=0pt] (3) at (4,0) {};
           \draw (1) -- (2);
            \node[anchor = east] at (0) {0};
           \node[anchor = south east] at (1,1) {1};
            \node[anchor = south] at (2) {2};
           \node[anchor = west] at (2,0) {\small{$4$}};
           \node at (1,0) {\small{$4$}};
           \node[anchor = north east] at (1,-1) {1};
           \node[anchor = north] at (1) {1};
           \node[anchor = west] at (3) {5};
           \node[anchor = south west] at (3,1) {6};
           \node[anchor = north west] at (3,-1) {6};
           \node at (3,0) {7};

            \node at (5,0) {$\leftarrow$};
            \node at (5,0.3) {$\varphi^*$};
            \draw[fill=lightgray] (0+6,0) -- (2+6,-2) -- (4+6,0) -- (2+6,2) -- cycle;
            \node[circle, fill=black, inner sep=1.5pt, outer sep=0pt] (a) at (0+6,0) {};
           \node[circle, fill=black, inner sep=1.5pt, outer sep=0pt] (b) at (2+6,-2) {};
           \node[circle, fill=black, inner sep=1.5pt, outer sep=0pt] (c) at (2+6,2) {};
           \node[circle, fill=black, inner sep=1.5pt, outer sep=0pt] (d) at (4+6,0) {};
           \draw (a) -- (b);
           \draw (a) -- (c);
           \draw (d) -- (b);
           \draw (d) -- (c);

           \node[anchor = east] at (a) {0};
           \node[anchor = south east] at (1+6,1) {1};
            \node[anchor = south] at (c) {2};
           \node[anchor = north east] at (1+6,-1) {1};
           \node[anchor = north] at (b) {1};
           \node[anchor = west] at (d) {5};
           \node[anchor = south west] at (3+6,1) {6};
           \node[anchor = north west] at (3+6,-1) {6};
           \node at (2+6,0) {7};
       \end{tikzpicture}
       \end{center}
   \end{Ex}
   \begin{Prop}
\label{Prop: pullbackpushforwardofinternalcollapses}       The pullback and push forward maps of internal collapses/extensions from \Cref{definition: pullback and push forward for internal collapses and extensions} are well-defined maps of spaces of discrete Morse functions.\par
       Moreover, the pullback map of an internal collapse, respectively the push forward map of an internal expansion, is an embedding of spaces of discrete Morse functions.
   \end{Prop}
   \begin{proof}
       For the pullback map induced by an internal collapse, taking the mean value in the construction ensures that the function values are decreasing along the newly created flow line, except for the pair of the internal collapse, on which $\varphi^*(f)$ attains the same value. This implies that the pair, which is being collapsed, is matched and $\varphi^*(f)$ is strictly monotone along all other (co-)face relations of $\sigma$ and $\tau$. It follows from $f$ being a discrete Morse function by assumption, that the definition of a discrete Morse function is satisfied at all other face relations as well. \par
       For the push forward map it is essential that we restrict to $\mathcal{M}(X)_{\vert x_\sigma\geq x_\tau}$, i.e. the regions of $\mathcal{M}(X)$ where $\sigma$ and $\tau$ are matched. Then it follows similarly as in the previous case, that $\varphi_*(f)$ fulfills the definition of a discrete Morse function at every face relation due to $f$ fulfilling it. \par
       For the last statement, we observe that pullbacks of internal collapses and push forwards of internal expansions leave values outside the pair $(\sigma,\tau)$ that is being collpased/expanded unchanged.
   \end{proof}
   \begin{Rem}\label{remark: space of dmfs on a simple homotopy type}
       The push forward and pullback maps induced by elementary and internal collapses/extensions allow us to relate spaces of discrete Morse functions on simple homotopy equivalent CW complexes with one-another. This way, it may make sense to develop the concept of a ``space of discrete Morse functions on a simple homotopy type" as a colimit of a diagram of spaces of discrete Morse functions on all representatives of a simple homotopy class with pullbacks and push forwards of all possible elementary and internal collapses and extensions between them. Due to the difficulties of such infinite diagrams and the combinatorial complexity of simple homotopy types, we leave such ideas for future work. 
   \end{Rem}
\subsection{The path metric on spaces of discrete Morse functions}
In the previous part, we did not explicitly mention how we topologize the spaces of discrete Morse functions. Instead, we rather implied that as a subspace of the finite dimensional vector space of discrete functions, we can use any topology, e.g.~the one induced by the Euclidean norm. \par
In \Cref{smoothcase} it will turn out that the Euclidean topology is useful for the relation to the smooth case. However, in applications towards topological data analysis it is often useful to be able to use several different specific metrics that induce the same topology. This provides the flexibility to used metrics adjusted to certain algorithmic approaches while other topological guarantees stay true due to preserving the induced topology. \par
In this subsection, we induce a specific notion of distance between discrete Morse functions on the same given regular CW complex $X$. This specific metric is meant to take the cost of moving through the hyperplanes of the braid arrangement $\mathcal{H}$ and the Morse arrangement $\mathcal{A}$ into account, and will be useful for the stability result \Cref{TheoInducedMergeTreeIsContinuous} for the induced merge tree. 
\begin{Rem}\label{Remark: PathsOnlyEnterOrLeaveIntersectionsOfHyperplanes}
	Let $X$ be a regular CW complex and let $\gamma \colon [0,1] \rightarrow \mathcal{M}(X)$ be a path in the space of discrete Morse functions. The path $\gamma$ cannot simultaneously enter an intersection of hyperplanes $\mathfrak{H}\coloneqq\bigcap\limits_j H_{\sigma_{j}}^{\tau_{j}}, H_{\sigma_j}^{\tau_j}\in \mathcal{H}$ and leave another one intersection of hyperplanes $\mathfrak{H}\coloneqq\bigcap\limits_{j'} H_{\sigma_{{j'}}}^{\tau_{{j'}}}, H_{\sigma_{j'}}^{\tau_{j'}}\in \mathcal{H}$ due to entering an intersection of hyperplanes is a closed condition and leaving intersections of hyperplanes is an open condition.  
\end{Rem}

\begin{Def}\label{DefinitionPathDistance}
    Let $X$ be a regular CW complex, let $\mathcal{M}(X)$ be the space of discrete Morse functions on $X$, and let $f,g\in \mathcal{M}(X)$.\par 
    Let $\Gamma \coloneqq \operatorname{Map}( [0,1]\rightarrow \mathcal{M}(X)) $ be the set of 
    paths $\gamma$ from $f$ to $g$ in the space of discrete Morse functions. 
    For any path $\gamma \in \Gamma$, let $\left(\mathfrak{H}_i= \bigcap\limits_j H_{\sigma_{i_j}}^{\tau_{i_j}}\right)$ be the sequence of intersections of hyperplanes that $\gamma$ crosses and let $\bar{t}\coloneqq(t_l)$ be an ordered sequence of times such that for each intersection of hyperplanes $\mathfrak{H}_i$ there is a $t_l$ such that $\gamma(t_l)\in \mathfrak{H}_i$ and for each face $F_i$ between two consecutive intersections of hyperplanes there is a $t_l$  such that $\gamma(t_l) \in F_i$. We call any such sequence $(t_l)$ and \D{admissible $\gamma$ decomposition} of $[0,1]$ and denote the collection of all admissible $\gamma$ decompositions by $T(\gamma)$. We define an ordered sequence $\bar{\gamma}\coloneqq(f\eqqcolon h_0, h_1,\dots, h_l, \dots,h_{k+1}\coloneqq g)$ of discrete Morse functions such that each $h_i$ for $1\leq i \leq k$  
    given by: $h_l \coloneqq \gamma(t_l)$ for $l\leq k$.
        
    We define the \D{path distance} on $\mathcal{M}(X)$ as follows: 
    $$
    d(f,g)\coloneqq \inf\limits_{\gamma \in \Gamma}\inf\limits_{\bar{t} \in T(\gamma)} \sum\limits_{l=0}^{k} \lvert \lvert h_l-h_{l+1}\rvert \rvert,
    $$
    where $\lvert\lvert \cdot\rvert\rvert$ denotes the Euclidean norm.
\end{Def}
\begin{Prop}
    The path distance is a metric.
\end{Prop}
\begin{proof}
    Reflexivity, positivity and symmetry follow directly from the definition. 
    The triangle inequality follows by concatenations of paths.
\end{proof}
\section{Relationship to other concepts in mathematics}
\subsection{Smooth Morse Theory and Cerf Theory}\label{smoothcase}
In this section we want to recall the description of the space of Morse functions on a smooth manifold given by Cerf in \cite{Cerf}. 
Moreover, we will explain the relationship between certain neighborhoods of Morse functions and the space of discrete Morse functions on their induced CW decompositions, inspired by \cite[Section 3.2]{Cerf}.\par
In order to get going, we start with a description of the space of smooth Morse functions on a given manifold $M$ and its relationship to Cerf's stratifications \cite{Cerf} of the space of smooth functions and the space of Morse functions. Throughout this section, we assume any smooth manifold $M$ to be finite-dimensional, compact, and with a fixed Riemannian metric.
\begin{Def}[\mbox{\cite[Definition 1,2, and 3]{Cerf}}]\label{Definition Cerf Strat}
Let $M$ be a smooth manifold. Denote by $C^\infty(M)$ the space of smooth functions from $M$ to $\mathbb{R}$ together with the $\mathcal{C}^\infty$ topology. Let $f\in C^\infty(M)$ be a smooth function, let $p\in M$ be a critical point of $f$, and let $a\in \mathbb{R}$ be a critical value of $f$. \par
The \D{codimension of $p\in M$} is defined as $\operatorname{codim}(\operatorname{\partial_p(f)} \subset C^\infty_{0,p}(M) )\coloneqq $ \\ $ \operatorname{dim}(C^\infty_{0,p}(M)/(\operatorname{\partial_p(f)}))$, where $C^\infty_{0,p}(M)$ denotes the ring\footnote{These rings do not necessarily contain a unit.} of germs of smooth functions vanishing at $p$ and  $(\operatorname{\partial_p(f)})$ denotes the ideal generated by the partial derivatives of $f$ at the point $p$.\par
The \D{codimension of $a$} is defined as $\operatorname{codim}(a)\coloneqq \lvert\{p\in f^{-1}(a) \vert p \text{ critical}\}\rvert-1$.\par
The \D{codimension of $f$} is defined as 
$$\operatorname{codim}(f)\coloneqq \sum\limits_{p\in M \text{ critical}}\operatorname{codim}(p)+\sum\limits_{a \in \mathbb{R} \text{ critical}} \operatorname{codim}(a).$$
Moreover, we define $\mathcal{F}\coloneqq C^\infty(M)$ and $\mathcal{F}^j\coloneqq \{ f \in \mathcal{F} \vert \operatorname{codim}(f) = j\}$. We call $\{\mathcal{F}^j\}_{j\in \mathbb{N}}$ the \D{natural stratification of $\mathcal{F}$}.
\end{Def}
We recall that for a function $f\in \mathcal{F}$ having a higher codimension than 0 means either that $f$ has degenerate critical points or multiple critical points of the same value. We denote the space of (smooth) Morse functions by $\mathcal{M}_{\text{smo}}(M)\subset \mathcal{F}$.\par
Next, we consider Cerf's comparison to the space of discrete functions.\par
We refer to \cite[Chapter 1, \S 3]{Sharko93} and \cite{Cerf} for background information on the natural stratification of the space of smooth functions $\mathcal{F}$. Before proceeding, we want to recall from \cite[Chapter 1, \S 3]{Sharko93} that $\mathcal{F}$ is a smoothly path connected smooth Fr\'echet manifold. From now on, all paths will be at least $C^1$ and path components will refer to $C^1$ path components.
\begin{Prop}[\mbox{\cite[Section 3.2: Definition 5, Lemma 1]{Cerf}}] \label{CerfLemma}
    Let $f\in \mathcal{M}_{\text{smo}}(M)$ be a smooth Morse function with $q$ critical points. Given a choice of an ordering of the critical points $c_1,\dots, c_q$, there are open neighborhoods $U_i$ with $c_i \in U_i\subset M$ for which there is an open neighborhood $\mathcal{V}$ with $f\in \mathcal{V}\subset \mathcal{F}$ such that all $\tilde{f}\in \mathcal{V}$ are Morse and have exactly one critical point $\tilde{c}_i\in U_i$ of the same index as $c_i$ for all $1\leq i \leq q$. Moreover, the ordering of the critical points of $f$ defines a topological submersion
    \begin{equation*}
    \begin{split}
    \eta \colon \mathcal{V} &\rightarrow \mathbb{R}^q \\
    \tilde{f} &\mapsto (\tilde{f}(\tilde{c}_1), \dots, \tilde{f}(\tilde{c}_q))  \end{split}
    \end{equation*}
    such that the restriction of the natural stratification of $\mathcal{F}$ to $\mathcal{V}$ is the preimage of the stratification induced by the hyperplane arrangement $\mathcal{H}_q$ (see \Cref{defbraidarrangement}). 
\end{Prop}
Our plan for this section is to identify the space $\mathbb{R}^q$ in \Cref{CerfLemma} with the space of discrete functions on the CW decomposition of $M$ induced by $f$, and use the stability of Morse--Smale functions to extend $\eta$ to path components. This allows for comparison maps from path components of the space of smooth Morse functions on $M$ to spaces of discrete Morse functions on certain cellular decompositions of $M$. For that, we recall the notion of Morse--Smale functions.
\begin{Def}[{\cite[Definition 6.1]{banyaga2004lectures}}]
    A Morse function $f\colon M \rightarrow \mathbb{R}$ on a finite dimensional smooth Riemannian manifold $(M,g)$ is said to satisfy the \D{Morse--Smale transversality condition} if and only if the stable manifold of $p$ and unstable manifold of $q$ with respect to $f$ intersect transversally\footnote{Recall that two submanifolds $N_1,N_2\subset M$ are said to intersect transversally if at each point of the intersection $p\in N_1 \cap N_2$ the tangent spaces of the submanifolds span the tangent space of the ambient manifold: $T_pN_1 + T_pN_2 = T_pM$.} for all pairs of critical points $p,q$ of $f$. A Morse function that satisfies the Morse--Smale transversality condition is called a \D{Morse--Smale function}.
\end{Def}
It is well known (see e.g. \cite[Theorem 6.34]{banyaga2004lectures}) that Morse functions only induce an explicit cellular decomposition $M_f$ of $M$, unique up to cell equivalence\footnote{Recall the \Cref{cell equivalence} of cell equivalences.}, if they fulfill the Morse--Smale property and are compatible with the given Riemannian metric. However, it is proved in \cite[Proposition 1.6]{FRANKS1979199} that the gradient of any Morse--Smale function can be perturbed to be compatible with the Riemannian metric without leaving its path component in the space of Morse--Smale vector fields, in particular while preserving topological conjugacy. Thus, we have a well-defined CW decomposition of $M$ induced by any Morse--Smale function. \par
Alternatively, one can associate gradient-like Morse--Smale vector fields to non-Morse--Smale functions instead (see \cite[Chapter 3 and Section 4.2]{Matsu}). Since we want to take the point of view of the space of Morse functions, we have fixed a Riemannian metric and consider the subspace $\mathcal{MS}(M) \subset \mathcal{M}_{\text{smo}}(M)$ of Morse--Smale functions on $M$.
\begin{Def}
    Let $Z$ be either the space of Morse functions $\mathcal{M}_{\text{smo}}(M)\subset \mathcal{F}$ or the space of Morse--Smale functions $\mathcal{MS}(M)\subset \mathcal{M}_{\text{smo}}(M)\subset \mathcal{F}$ and let $f\in Z$. We denote by $\mathcal{N}(f)\subset Z$ the path component of $f$ in $Z$\footnote{It will be clear from the context which of the two possibilities for $Z$ will be used.}. 
\end{Def}
It is a classical result that the property of being a Morse function is stable, i.e.\ for any Morse function $f$, there is an open neighborhood $U(f)\subset C^\infty(M)$ such that $U(f)\subset \mathcal{M}_{\text{smo}}(M)$ and for every path $f_\_\colon [0,1] \rightarrow U(f)$ there is an $\epsilon >0$ such that $f_t$ is Morse for every $t<\epsilon$ (see \cite[Corollary 5.24]{banyaga2004lectures}). This classical result is also used in the proof of $\Cref{CerfLemma}$. Moreover, in \cite[\S 1]{FRANKS1979199} the stability of Morse--Smale functions is stated in an even stronger way, which we use for the following theorem: 
\begin{Def}
Let $\phi,\phi'$ be flows on a closed orientable Riemannian manifold $M$. $\phi$ and $\phi'$ are called \D{topologically conjugate} if there exists a homeomorphism $M\rightarrow M$ that carries orbits of $\phi$ to orbits of $\phi'$ and preserves their orientation. \par
\end{Def}
\begin{Theo}[{\cite[Theorem 1.12]{PALIS1969385}, reformulation from \cite[\S 1]{FRANKS1979199}}]
Let $\phi$ be a Morse--Smale flow on a closed orientable Riemannian manifold $M$ and let $\phi'$ be a sufficiently small $C^1$ perturbation of $\phi$, then $\phi$ and $\phi'$ are topologically conjugate.    
\end{Theo}
For the proof of this theorem, we refer to \cite[Theorem 1.12]{PALIS1969385}.
\begin{Theo}\label{maintheosmooth}
    Let $f,g$ be Morse--Smale functions on $M$ which are contained in the same path component, i.e.\ $g\in \mathcal{N}(f)\subset \mathcal{MS}(M)$. Then there is a bijection $\psi \colon \operatorname{Cr}(f)\xrightarrow[]{\cong} \operatorname{Cr}(g)$ such that there is an orbit of $\phi_f$ connecting two critical points $c_1,c_2 \in \operatorname{Cr}(f)$ if and only if there is an orbit of $\phi_g$ connecting $\psi(c_1),\psi(c_2)\in \operatorname{Cr}(g)$.
\end{Theo}
\begin{proof}
    This result is a consequence of the structural stability of Morse--Smale functions, respectively Morse--Smale vector fields, respectively Morse--Smale flows. 
    Let $f_t$ be a path from $f$ to $g$ in $\mathcal{MS}(M)$. We choose an open cover on $I$ that consists of such open neighborhoods $U(t_i)$ of time steps $t_i$ that the restriction of $f_t$ to $U(t_i)$ is a small enough $C^1$ perturbation so that stability in the sense of \cite{FRANKS1979199} holds. Due to $I$ being compact, there is a finite subcover of neighborhoods where the restriction of $f_t$ induces topologically conjugate flows. Thus, $f_t$ induces topologically conjugate flows for all $t$. Hence, the statement follows. 
    \end{proof}

\begin{Cor}\label{Decompositionwelldefined}
    For any compact finite-dimensional Riemannian manifold $M$, all Morse--Smale functions which are in the same path component of $\mathcal{MS}(M)$ induce cell-equivalent CW decompositions $M_f$ of $M$. That is, there is an up to cell equivalence well-defined CW decomposition $M_\mathcal{N}$ associated to any path component $\mathcal{N}$ of the space of Morse--Smale functions on $M$. 
\end{Cor}
\begin{Cor}\label{CorCerf}
    It follows from \Cref{maintheosmooth} that Cerf's map from \Cref{CerfLemma} extends to a map $\eta \colon \mathcal{N}\rightarrow \mathcal{M}(M_{\mathcal{N}})$, where $\mathcal{M}(M_{\mathcal{N}})$ denotes the CW decomposition of $M$ induced by any Morse function in the path component $\mathcal{N}$ of $\mathcal{M}(M)$. Moreover, Cerf's proof of \Cref{CerfLemma} shows that this instance of $\eta$ is a topological submersion and compatible with the respective stratifications, too.
\end{Cor}
Building on \Cref{Decompositionwelldefined}, we can extend the map $\eta$ from \Cref{CerfLemma} to the path component $\mathcal{N}(f)$ of a given Morse--Smale function $f$. By \Cref{CerfLemma}, it still follows that $\eta \colon \mathcal{N}(f)\rightarrow \mathcal{M}(M_{\mathcal{N}(f)}) \coloneqq \mathcal{M}(M_f)$ is a topological submersion which is compatible with the respective stratifications.\par
Using Cerf's stratification of the space of Morse functions, especially the connection between different path components given by cancellations and creations of pairs of critical points, we want to further extend $\eta$ to certain systems of path components:
\begin{Def}\label{cancellationorder}
    Consider the set $\pi_0(\mathcal{MS}(M))$ of path components of the space of Morse--Smale functions on $M$. Let $\mathcal{N}_1,\mathcal{N}_2\in \pi_0(\mathcal{MS}(M))$ be path components. Then we define a partial order $\leq$ on $\pi_0(\mathcal{MS}(M))$ generated by $\mathcal{N}_1 < \mathcal{N}_2$ if for any $f\in \mathcal{N}_2$ there is a $g\in \mathcal{N}_1$ such that $g$ arises from $f$ by a cancellation of a pair of critical points. We call $\leq$ the \D{cancellation order} on $\pi_0(\mathcal{MS}(M))$.
\end{Def}
\begin{Rem}
    It is shown in \cite[Section 3]{Cerf} that if two generic Morse functions, i.e.~elements of $\mathcal{F}^0$, are related to each other by a cancellation of one pair of critical points, then all paths in $\mathcal{F}$ between them are homotopic in $\mathcal{F}$ to one which traverses the first stratum  $\mathcal{F}^1$ in one isolated point and is otherwise entirely contained in $\mathcal{F}^0$. Hence, we can think about two path components $\mathcal{N}_1,\mathcal{N}_2$ such that $\mathcal{N}_1 <\mathcal{N}_2$ as being neighbors, only separated by a codimension one stratum which corresponds to a hyperplane in $\mathcal{M}(\mathcal{N}_2)$ via $\eta$. 
\end{Rem}
\begin{Prop}\label{cancellationsimplecollapse}
    Let $f,g\in \mathcal{MS}(M)$ be two Morse--Smale functions such that $f$ arises from $g$ by a cancellation of one pair of critical points. Then $M_f$ and $M_g$ are simply homotopy equivalent by a simple collapse $M_g\rightarrow M_f$/ a simple extension $M_f\rightarrow M_g$ of the two cells which correspond to the canceled critical points at the corresponding sublevel complex. 
\end{Prop}
\begin{proof}
    This statement is basically a rephrased version of \cite[Theorem 3.34 (canceling handles: rephrased)]{Matsu}. In \cite[Theorem 3.34]{Matsu}, building on techniques from \cite[Proof of Theorem 5.4]{Milnor+1965} the authors argue that a cancellation of a pair of critical points induces a pair of handles in the induced handle body decomposition such that the handlebody decomposition induced by $f$ is diffeomorphic to the one induced by $g$. The statement then follows in a straightforward manner from the construction of the CW decomposition associated to a handlebody decomposition \cite[Theorem 4.18]{Matsu}.
\end{proof}
\begin{Cor}\label{Cor: MapofSpaceofDMFINducedByCancellation}
    Let $\mathcal{N}_1<\mathcal{N}_2\in \pi_0(\mathcal{MS}(M))$ be path components, which are related by cancellations of critical points. Then $\mathcal{A}(M_{\mathcal{N}_1})$ canonically embeds into $\mathcal{A}(M_{\mathcal{N}_2})$ and we have an embedding $\mathcal{M}(M_{\mathcal{N}_1})\rightarrow\mathcal{M}(M_{\mathcal{N}_2})$ by \Cref{Prop: pullbackpushforwardofinternalcollapses}.
\end{Cor}
\begin{Rem}
    Although the map $\eta$ from \Cref{CorCerf} is only defined on the path components of the space of Morse--Smale functions, \Cref{Cor: MapofSpaceofDMFINducedByCancellation} allows us to relate the spaces of discrete Morse functions on CW decompositions induced from different path components to each other. Provided one finds a good notion of space of discrete Morse functions on a simple homotopy type, as we have conjectured in \Cref{remark: space of dmfs on a simple homotopy type}, one could ``globalize" the map $\eta$ and extend it to a topological submersion $\eta \colon \mathcal{MS}(M) \rightarrow \mathcal{M}(\mathcal{SH}(M))$, where $\mathcal{SH}(M)$ denotes the simple homotopy type of all CW decompositions of $M$.  
\end{Rem}
\subsubsection{Discrete Cerf theory}
Since one of the central works that inspired this one is Cerf's structural description of the geometric interpretation of cancellation moves in the space of Morse functions, see \cite{Cerf}, we want to provide a brief translation of this theory to the much easier discrete case.\par
\begin{figure}[H]
    \centering
    \begin{tabularx}{\textwidth}{X|X}
       Smooth case & Discrete case\\ \hline
        -stratification by codimension & -stratification by braid arrangement\\
       -function has a critical point of codim 1& -function lies inside a hyperplane of the Morse arrangement\\
       -function has a critical point of higher codimension & -function lies inside the intersection of hyperplanes whose cover relations share a cell \\
       -universal unfolding of a critical point & -basis of the orthogonal complement of the (intersection of) hyperplane(s)\\
       -cancellation of a pair of critical points &-cancellation of a pair of critical cells\\
       -creation of a pair of critical points &-removing a matched pair from the Morse matching (potential need for an internal extension) \\
       -function has a critical value of higher codimension & -function lies in the intersection of hyperplanes that belong to the braid arrangement but not to the Morse arrangement\\
       -moving critical points past each other &-crossing a hyperplane that belongs to the braid arrangement but not to the Morse arrangement\\
    \end{tabularx}
    \caption{Comparison between smooth and discrete Cerf theory}
    \label{fig:placeholder}
\end{figure}
The well-known Cerf diagrams apply the same way in the discrete case as they do in the smooth case. Moreover, the Cerf diagrams can in the discrete case quite explicitly illustrate the corresponding path of discrete functions, seen inside the correct orthogonal projection of the space of discrete functions. \par
The central difference between the two theories is, that the creation of pairs of critical points is always available in the smooth case, whereas in the discrete case the CW complex at hand may be lacking the necessary cells. In order to fix this, one can always perform an internal extension, but this technically changes the space of discrete functions under consideration. In order to truly fix this, one would need to make sense of the notion of space of discrete Morse functions on a simple homotopy type, as we conjectured in \Cref{remark: space of dmfs on a simple homotopy type}.
\subsection{The space of discrete Morse matchings}
In the literature, there is already work on a notion of complex of discrete Morse functions. In fact, this name has been used with a slight abuse of terminology since these works refer to a complex whose simplices correspond to Morse matchings rather than Morse functions. In order to make this difference clear, we use the different name ``complex of discrete Morse matchings" to refer to this concept.
\begin{Def}[\mbox{\cite[Section 2]{ChariJo}}]\label{Def:spaceOfMorseMatchings}
Let $X$ be a simplicial complex. The \D{complex of discrete Morse matchings} $\mathfrak{M}(X)$ has as vertices cover relations in $F(X)$ and simplices sets of cover relations in $F(X)$ that correspond to acyclic matchings on $D(X)$.  
\end{Def}
It is straightforward how this notion should be generalized to arbitrary CW complexes:
\begin{Def}\label{DefcomplexofDMMs}
    Let $X$ be a CW complex. The \D{complex of discrete Morse matchings} $\mathfrak{M}(X)$ on $X$ has vertices cover relations in $F(X)$ that correspond to regular faces and simplices sets of vertices that correspond to acyclic matchings on $D(X)$. 
\end{Def}
Since the geometric structure of the simplices, i.e.~their points, do not have any obvious meaning for discrete Morse matchings, we will not distinguish between $\mathfrak{M}(X)$ and its face poset $F(\mathfrak{M}(X))$ in our notation. 
\begin{Prop}\label{PropSpaceOfMorseMatchings}
    Let $X$ be a CW complex and let $\mathcal{A}$ be the Morse arrangement on $\mathbb{R}^X$. Then there is a canonical embedding of posets $\mathfrak{M}(X)\subset \mathcal{L}(\mathcal{A})$.
\end{Prop}
\begin{proof}
    The minimal elements of both posets are in canonical bijection because in both cases they are defined by the regular cover relations in $F(X)$. Then collections of cover relations that form an acyclic matching in $\mathfrak{M}(X)$ are mapped to intersections of hyperplanes in $\mathcal{L}(\mathcal{A})$ that correspond to the same acyclic matching. 
\end{proof}
\subsection{The space of merge trees}
    In order to give a model for the space of merge trees, we introduce the notions of edit moves and edit morphisms to relate different combinatorial merge trees to each other.
    \begin{Def}\label{editmoves}
        An \D{edit move} on a poset of merge tree type $P$ is one of the following operations:
        \begin{enumerate}
            \item[The identity:]\mbox{}\\
            The poset $P$ is left unchanged. We denote this edit move by $\operatorname{id}$.
            \item[Adding a leaf:]\mbox{}\\
            A new minimal element $x$ is added in one of two ways: either, we just introduce the relation $x\leq y$ for a pre-existing non-minimal element $y\in P$, or we add a new non-minimal element $y'$ between two pre-existing elements $z_1\leq y' \leq z_2$ and add the relation $x\leq y'$. As an exceptional case, if $P$ is the one element poset $\{r\}$, in particular, there are no inner nodes, we allow to add two leaves $x<r,y<r$ simultaneously. We denote this edit move by $\operatorname{Add}(x,y)$.
            \item[Removing a leaf:]\mbox{}\\
            A pre-existing minimal element $x\in P$ is removed. If $x$'s parent node $y$ has only the one child node $x$, we remove $y$ in addition to $x$. We denote this by $\operatorname{Rem}(x,y)$.
            \item[Splitting an inner node]\footnote{On combinatorial merge trees, this type of edit move provides multiple different options of edit moves depending on how the child nodes are precisely being distributed among $z_1$ and $z_2$.}\mbox{}\\
            If a node $z\in P$ has at least 3 child nodes, then $z$ can be split into two nodes $z_1\leq z_2$ such that there is no $x\in P$ with $z_1\leq x\leq z_2$. The previous child nodes of $z$ are then distributed among $z_1$ and $z_2$ such that $z_1$ gets at least two child nodes and $z_2$ gets at least one child node other than $z_1$. We denote this edit move by $\operatorname{Spl}(z,z_1,z_2)$.
            \item[Merging two adjacent inner nodes:]\mbox{}\\ 
            If there are two non-minimal elements $z_1\leq z_2\in P$ such that $z_2$ covers $z_1$, i.e.~$z_1<z_2$ and there exists no $x: z_1\leq x\leq z_2$, we construct the quotient $P/_{z_1=z_2}$. This creates an element $z$ instead of $z_1$ and $z_2$ such that all previous child nodes of $z_1$ and $z_2$ become child nodes of $z$. We denote this edit move by $\operatorname{Mer}(z_1,z_2,z)$.            
        \end{enumerate}
        An \D{edit morphism} $\eta \colon P \rightarrow P'$ between posets of merge tree type is a composition of edit moves composed with an automorphism of $P'$. An edit morphism is called \D{elementary} if it performs only one of the edit moves. An edit morphism is called \D{trivial} if it has the identity as its only corresponding edit move. We some times abbreviate the notion \D{elementary edit morphism} by \D{EEM}. By slight abuse of notation, we denote elementary edit morphisms by their corresponding edit move. 
    \end{Def}
    \begin{Rem}\label{editmorphisminducedmap}
        Let $P,P'$ be posets of merge tree type such that $P'$ arises from $P$ by either an addition of a leaf or a splitting of an inner node. Then any corresponding edit morphism induces a poset morphism $\eta \colon P \rightarrow P'$ given by the inclusion of the preexisting nodes and mapping nodes $x$ that are being split to maximal one $z_2$ of the new nodes.\par
        Additionally, we want to remark that for each finite poset of merge tree type $P$, there is only a finite number of elementary edit moves that can be applied to $P$. Those are:
        \begin{enumerate}
            \item Adding one leaf to any inner node or edge,
            \item Removing at most one leaf per leaf (in general not possible for all leaves.)
        \end{enumerate}
    \end{Rem}


We identify spaces of discrete functions, as introduced in \Cref{definitionspaceofdiscretefunctions}, with sets of set maps $\operatorname{Map}(P,\mathbb{R})$ for posets $P$ of CW type. Moreover, we identify sets of merge trees with the same underlying combinatorial merge tree $P$ with sets of order-preserving maps $\operatorname{Map}_\leq(P,\mathbb{R})$ for posets $P$ of merge tree type. We use the following model\footnote{By model we mean a set that has one specific poset of merge tree typ for every isomorphism class of posets of merge tree type.} for the class of posets of merge tree type:
\begin{Def}
    We construct the \D{set of combinatorial merge trees} $Mer_{comb}$ as follows:
    Let $n\in \mathbb{N}_{\geq1}$. We first define a set of representatives for posets of merge tree type with $n$ nodes, which we denote by $Mer_{comb}^n$: we consider the set of nodes $\{1,2,\dots,n\}$ and fix a collection of poset structures on $\{1,2,\dots,n\}$ such that there is one representative for each isomorphism type of posets of merge tree type with $n$ nodes. This is possible because the collection of partial orders on $\{1,\dots,n\}$ is a subset of the set of relations $\{1,\dots,n\}\times\{1,\dots,n\}$, so in particular there is only a finite amount of those. We define $$Mer_{comb}\coloneqq \coprod\limits_{n\in \mathbb{N}_{\geq1}} Mer_{comb}^n.$$
\end{Def}
We do \textbf{not} claim that the choice of representatives in $Mer_{comb}$ is canonical in any way - \textbf{it is not}. However, this is fine for our purposes. From now on, we will not address this choice of representatives any further: we will always assume that composition with an additional poset isomorphism from the target of any morphism under consideration to one of the chosen representatives takes care of all such issues.  
    
\begin{Def}\label{DefSetsOfMergeTrees}
    We define the \D{set of merge trees} $Mer$ as the disjoint union $$Mer\coloneqq \coprod\limits_{P \in Mer_{comb}} \operatorname{Map}_\leq(P,\mathbb{R}).$$
    For any merge tree $\theta \in Mer$ we denote the underlying combinatorial merge tree, i.e.~the source of $\theta$, by $s(\theta)$. \par
    Moreover, we define the \D{set of strict merge trees} $Mer_{st}$ by  $$Mer_{st}\coloneqq \coprod\limits_{P Mer_{comb}} \operatorname{Map}_<(P,\mathbb{R}),$$ i.e.~we here demand that the height functions on the combinatorial merge trees are strictly monotone. \par
    Furthermore, we define the \D{set of well-branched merge trees} $Mer_{wb} \subset Mer_{st}$ by $$Mer_{wb}\coloneqq \{\theta \in Mer_{st} \vert \ \forall y \in s(\theta) \text{ the restriction } \theta_{\lvert s(\theta)_{\leq y}} \text{ has a minimum}\},$$
     where $s(\theta)_{\leq y}$ denotes the rooted subtree of $s(\theta)$ with root $y$.
\end{Def}
\begin{Rem}
    In the definition above we slightly abuse notation by identifying a merge tree with its height function. Since any height function $\theta$ implicitly also contains the information of the underlying combinatorial merge tree as the domain of $\theta$, this should hopefully not lead to any confusion.
\end{Rem}
We consider how edit morphisms can be used to shift height functions from one combinatorial merge tree to another:
\begin{Def}
We say that a collection of edit morphisms $\{\eta_i\}_{i \in I}$ is/goes \D{in the same direction} if there is a way to compose all of the $\eta_i$ to one morphism and all $\{\eta_i\}_{i \in I}$ are either additions of leaves and splittings of inner nodes, or removals of leaves and merges of inner nodes. \par
Let $\eta$ be a composition of edit morphisms in the same direction that consists of additions of leaves and splittings of inner nodes. We denote by $\operatorname{L}(\eta)$ the collection of pairs $(x,y)$ such that $x$ is a leave that is being added by $\eta$ to its parent node $y$. Moreover, we denote by $\operatorname{S}(\eta)$ the collection of pairs of nodes $(x,x_i)$ such that $x_i$ was created by $\eta$ splitting $x_i$ off of $x$. 
\end{Def}
In the case of merge trees whose underlying combinatorial merge trees are related by compositions of elementary edit morphisms in the same direction, there is a convenient way of shifting height function values between combinatorial merge trees:
\begin{Def}\label{DefinitionShiftOfHeightFunctionViaEEMSInTheSameDirection}
Let $\theta \in \operatorname{Map}_\leq(P,\mathbb{R})$
be a merge tree $P$ and let $P'$ be a poset of merge tree type such there is a composition of EEMs in the same direction $\eta\colon P \rightarrow P'$. \par
Then by \Cref{editmorphisminducedmap} we have an induced injective morphism $\eta \colon P \rightarrow P'$ if $\eta$ is a composition of adding leaves and splitting inner nodes. 
    We extend the inclusion $\eta \colon P \rightarrow P'$ associated to any composition of adding leaves and splitting inner nodes to a map $\eta_* \colon \operatorname{Map}_\leq(P,\mathbb{R}) \times \operatorname{Map}_\leq(P',\mathbb{R}) \rightarrow \operatorname{Map}_\leq(P',\mathbb{R})$ by  
    \begin{align}\label{FormulaInducedMapEta}
    \eta_*(\theta,\theta')_x\coloneqq \begin{cases}
        \theta_{\eta^{-1}(x)} & \text{ if } x \in \eta(P),\\
        \theta_{\eta^{-1}(y)} & \text{ if } x \notin \eta(P), \ x \text{ a leaf, } y \in \eta(P) \text{ is }  x \text{'s parent  in } P'\\
        \theta'_{x} & \text{ if } x \notin \eta(P) \text{ was added as a new parent node for a new leaf}\\
        \theta'_y & \text{ if } x \notin \eta(P), \ x \text{ a leaf, } y \notin \eta(P) \text{ is }  x \text{'s parent  in } P'\\
        \theta_y & \text{ if } x \notin \eta(P), \ x \text{ an inner node, } y  \text{ the node } x \text{ was split from}
    \end{cases}. \end{align}
\end{Def}
\begin{Rem}
    If $\eta$ is a composition of removing leaves and merging inner nodes, we have an induced inclusion $\eta\colon P' \rightarrow P$ in the other direction by the same construction for the inverse composition of EEMs in the same direction.
\end{Rem}
\begin{Def}\label{inducedmergetreemap}
    Let $X$ be a regular CW complex. We define a map $M \colon \mathcal{M}(X) \rightarrow Mer$ as follows: Let $f\in \mathcal{M}(X)$ be a discrete Morse function and let $\operatorname{Cr}_i(f)$ the set of critical cells of dimension $i$ of $f$. We create a poset of merge tree type $\tilde{P}$ as follows: for each critical vertex $\sigma \in \operatorname{Cr}_0(f)$ we add an element $\sigma$ to $\tilde{P}$. For each time when connected components $\sigma_1,\dots, \sigma_k$ of sublevel complexes merge within the filtration, we add an element $\tau$ to $\tilde{P}$ and add the relations $\sigma_i\leq \tau$. The elements $\tau$ correspond to certain critical edges $\tau \in \operatorname{Cr}_1(f)$. Moreover, we identify $\tilde{P}$ with the unique representative $P\in Mer_{comb}$ of the isomorphism type of $\tilde{P}$. Then we define $M(f)(\sigma) \coloneqq f_\sigma$, i.e.~the elements of $P$ get their function values under $M(f)$ from the critical values of the corresponding critical cells. 
\end{Def}
\begin{Ex}\label{ExampleInducedMergeTree}
We consider an example of a discrete Morse function on a regular CW complex and its induced merge tree:
    \begin{center}
    \begin{tikzpicture}
    \draw[fill=lightgray] (5,-1) -- (6.5,-2) -- (6.5,0.25) -- cycle;
    \draw[color=pink] (5,-1) -- (6.5,-2);
    \draw[color=lime] (6.5,-2) -- (6.5,0.25);
	\shade[ball color = gray!40, opacity = 0.4] (3.5,-1) circle (1.5cm);
	\draw (3.5,-1) circle (1.5cm);
	\draw[color=orange] (2,-1) arc (180:360:1.5 and 0.6);
	\draw[dashed,color=orange] (5,-1) arc (0:180:1.5 and 0.6);
    \draw (2,-1) node {$\bullet$};
	\draw (5,-1) node[color=blue] {$\bullet$};
    \node[circle, fill=green, inner sep=1.5pt, outer sep=0pt] (0) at (6.5,-2) {};
    \node[circle, fill=purple, inner sep=1.5pt, outer sep=0pt] (1) at (6.5,0.25) {};
	\draw (2,-1) node[anchor=east] {0};
    \draw (3.5,-1) node[anchor=north] {2};
    \draw (3.5,-1) node[anchor=south] {2};
    \draw (3.5,0.5) node[anchor=south] {3};
    \draw (3.5,-2.5) node[anchor=north] {4};
    \draw (4.9,-1.1) node[anchor=north west] {1};
    \draw (5.75,-0.375) node[anchor=south] {5};
    \draw (5.75,-1.625) node[anchor=north] {4};
    \draw (5.8,-1) node[] {5};
    \draw (6.5,0.25) node[anchor=south] {0};
    \draw (6.5,-2) node[anchor=north] {1};
    \draw (6.5,-1) node[anchor=west] {3};

    \node[circle, fill=pink, inner sep=1.5pt, outer sep=0pt] (a) at (10,0.25) {};
    \draw (10,0.25) node[anchor=south] {4};
    \node[circle, fill=orange, inner sep=1.5pt, outer sep=0pt] (b) at (9,-0.75) {};
    \draw (9,-0.75) node[anchor=east] {2};
    \node[circle, fill=lime, inner sep=1.5pt, outer sep=0pt] (c) at (11,-0.75) {};
    \draw (11,-0.75) node[anchor=west] {3};
    \node[circle, fill=black, inner sep=1.5pt, outer sep=0pt] (d) at (8.5,-1.75) {};
    \draw (8.5,-1.75) node[anchor=north] {0};
    \node[circle, fill=blue, inner sep=1.5pt, outer sep=0pt] (e) at (9.5,-1.75) {};
    \draw (9.5,-1.75) node[anchor=north] {1};
    \node[circle, fill=green, inner sep=1.5pt, outer sep=0pt] (f) at (10.5,-1.75) {};
    \draw (10.5,-1.75) node[anchor=north] {1};
    \node[circle, fill=purple, inner sep=1.5pt, outer sep=0pt] (g) at (11.5,-1.75) {};
    \draw (11.5,-1.75) node[anchor=north] {0};
    \draw (b) -- (a);
    \draw (c) -- (a);
    \draw (b) -- (d);
    \draw (b) -- (e);
    \draw (c) -- (f);
    \draw (c) -- (g);
    \end{tikzpicture}
    \end{center}
    The colors indicate which critical cells correspond to which nodes of the merge tree. Here we have drawn the merge tree as its Hasse diagram.
\end{Ex}
\begin{Prop}
    The image of $\mathcal{MB}_w(X)$ is contained in $Mer_{wb}$.
\end{Prop}
\begin{proof}
    Since weak Morse--Benedetti functions cannot obtain the same value on different critical cells, the even stronger statement holds that $M(\mathcal{MB}_w(X))$ only contains injective maps for any regular CW complex $X$.
\end{proof}
We proceed by considering how crossing hyperplanes of $\mathcal{A}$ and of $\mathcal{H}$ relates to edit moves between the induced combinatorial merge trees.
\begin{Lem}\label{Lemmacancelcrossing}
    Let $X$ be a regular CW complex and let $H_\sigma^\tau$ be a hyperplane of $\mathcal{A}(X)$. If $\operatorname{dim}(\sigma)=0$, then passing from a Morse region $R_+$ on the positive side of $H_\sigma^\tau$ to a Morse region $R_-$ on the negative side without crossing any other hyperplane of $\mathcal{H}$ corresponds to the removal of a leaf. Moreover, passing through $H_\sigma^\tau$ the other way around corresponds to adding a leaf.
\end{Lem}
\begin{Rem}
    We observe that a crossing as mentioned in \Cref{Lemmacancelcrossing} does not exist for all Morse regions on the positive side of $H_\sigma^\tau$. Even if it is possible to cancel $\sigma$ and $\tau$ without violating the Morse condition, it might happen that one has to cross other regions of $\mathcal{H}$ or even other Morse regions along the way. 
\end{Rem}
\begin{proof}[Proof of \Cref{Lemmacancelcrossing}]
    Since, by assumption, both $R_+$ and $R_-$ are Morse regions and it is possible to pass from $R_+$ to $R_-$ without crossing any hyperplane of $\mathcal{H}$ other than $H_\sigma^\tau$, this means that for any Morse function $f\in R_+$, we have $f(\sigma)<f(\tau)$ and there is no cell with a function value between them. Since $\operatorname{dim}(\sigma)=0$, at filtration level $f(\sigma)$ a new connected component arises, which is merged with another connected component at level $f(\tau)$ via $\tau$. Because passing to $R_-$ corresponds to matching $\sigma$ and $\tau$, moving to $R_-$ removes the leaf of $M(f)$ that corresponds to $\sigma$. It is immediate that passing from $R_-$ to $R_+$ adds a new leaf to the induced merge tree.
\end{proof}
\begin{Rem}\label{Remcancelpath}
    It is well known, e.g. \cite[Theorem 11.1]{Forman98}, that whenever a Morse function $f\in \mathcal{M}(X)$ induces a unique gradient path between a critical $d$ cell $\sigma$ and a critical $d+1$ cell $\tau$ such that all face relations along the gradient path are regular, the gradient path can be inverted, matching $\sigma$ and $\tau$ as a result. If $\sigma$ and $\tau$ have consecutive function values, then $\sigma$ is a face of $\tau$ and we are in the situation of \Cref{Lemmacancelcrossing}. Otherwise, there are other cells with function values between $f(\sigma)$ and $f(\tau)$ and we have to cross hyperplanes of $\mathcal{A}$ corresponding to all face relations along said unique gradient path between $\sigma$ and $\tau$ as well as hyperplanes of $\mathcal{H}\setminus\mathcal{A}$ corresponding to cells that are not contained in the gradient path with function values between $f(\sigma)$ and $f(\tau)$. It should be remarked that during this procedure, we cannot cross all hyperplanes of $\mathcal{A}$ corresponding to the face relations along the unique gradient path simultaneously, if the path is of length larger than 1, because otherwise we would leave the space of discrete Morse functions.  
\end{Rem} 
\begin{Lem}\label{LemmaMoveCrossing}
    Let $X$ be a regular CW complex, let $\sigma,\tau\in X$ be two 1 cells, and let $H_\sigma^\tau$ be the corresponding hyperplane in $\mathcal{H}$. Let $f\in \mathcal{M}$ be a discrete Morse function inside a region $R_1$ of $\mathcal{H}$ such that $\sigma$ and $\tau$ are critical and merge connected components, i.e.~correspond to inner nodes of $M(f)$. Assume there is a region $R_2$ adjacent to $R_1$ via $H_\sigma^\tau \in \mathcal{H}$, i.e.~there is a path from $R_1$ to $R_2$ that only crosses $H_\sigma^\tau$ and no other hyperplane of $\mathcal{H}$. Then the following statements hold:
    \begin{enumerate}
        \item If $\sigma$ and $\tau$ correspond to adjacent inner nodes of $M(f)$, then moving through $H_\sigma^\tau$ first merges the corresponding nodes and then splits them the other way around.
        \item If the nodes corresponding to $\sigma$ and $\tau$ are not comparable in $M(f)$, passing through $H_\sigma^\tau$ corresponds to moving their associated values past each other without changing the underlying combinatorial merge tree. 
    \end{enumerate}
    
\end{Lem}
    \begin{Rem}
        In the context of \Cref{LemmaMoveCrossing}, if the nodes that correspond to $\sigma$ and $\tau$ are comparable but not adjacent in $M(f)$, then it is not possible to cross $H_\sigma^\tau$ without crossing other hyperplanes of $\mathcal{H}$ because then there are other nodes with values in between due to the monotonicity of values associated to nodes of $M(f)$.
    \end{Rem}
    \begin{proof}[Proof of \Cref{LemmaMoveCrossing}]
        \begin{enumerate}
            \item If $\sigma$ and $\tau$ correspond to adjacent inner nodes of the induced merge tree, then $\sigma$ and $\tau$ subsequently merge three previously disconnected connected components into one. Crossing $H_\sigma^\tau$ means that the components are merged simultaneously while inside $H_\sigma^\tau$ and the merge order is reversed after exiting $H_\sigma^\tau$ on the other side. This corresponds by definition to first merge the nodes and then splitting them the other way around.
            \item If the nodes corresponding to $\sigma$ and $\tau$ are not comparable then $\sigma$ and $\tau$ belong to different connected components which are, if at all, merged at higher level via a third 1-cell. Hence, crossing $H_\sigma^\tau$ only changes the levels at which $\sigma$ and $\tau$ merge their respective components but not the underlying combinatorial merge tree. 
        \end{enumerate}
    \end{proof}
We define three hyperplane arrangements on the space of merge trees, which are related to the Morse arrangement and the braid arrangement on the space of Morse functions on any given CW complex:
\begin{Def}
    Let $P$ be a poset of merge tree type. We define the \D{leaf arrangement} $\mathcal{A}(P)$ in $\operatorname{Map}_\leq(P,\mathbb{R})$ by $\mathcal{A}(P)\coloneqq \{H_x^y \vert x \text{ is a leaf and } y \text{ is the parent node of }x \}$.\par
    We define the \D{order arrangement} $\mathcal{O}(P)$ by $\mathcal{O}(P)\coloneqq \{H_x^y \vert x\leq y\}$
    .\par
    Moreover, we consider the \D{braid arrangement} $\mathcal{H}(P)$ in $\operatorname{Map}_\leq(P,\mathbb{R})$, which is given by $\mathcal{H}(P) \coloneqq \{H_x^y \vert x,y \in P \}$. 
\end{Def}

\begin{Prop}\label{PropGradientPahtsAndHyperplanes}
    Let $P$ be a poset of merge tree type and let $f\colon X \rightarrow \mathbb{R}$ a discrete Morse function on a CW complex $X$ such that $P$ is isomorphic to the underlying combinatorial merge tree of $M(f)$. Then the following holds:
    \begin{enumerate}
    \item  The hyperplanes of $\mathcal{A}(P)$ correspond to families of maximal gradient paths between critical 0 and 1 cells induced by $f$. If $f$ is weakly Morse--Benedetti, then this correspondence is bijective. 
    \item The hyperplanes of $\mathcal{O}(P)$ correspond to a collection of zig-zags of maximal gradient paths between  cells of dimension $\leq 1$ induced by $f$ that begin and end at maximal critical cells of sublevel complexes without leaving the sublevel complex corresponding to the higher of the two cells.
    \item The hyperplanes of $\mathcal{H}(P)$ correspond to a collection of zig-zags of maximal gradient paths between cells of dimension $\leq 1$ induced by $f$.
    \end{enumerate}
    In particular $\mathcal{A}(P)\subset \mathcal{O}(P)\subset \mathcal{H}(P)$ holds.
\end{Prop}
\begin{proof}
    \begin{enumerate}
        \item By \Cref{inducedmergetreemap}, the leaves of $P$ correspond to critical 0 cells of $f$. Moreover, the parent relation $x\prec y$ of any leaf $x$ corresponds to the collection of critical 1 cells that merge the connected component that corresponds to $x$, say $X_{M(f)_{y}-\varepsilon}[\sigma]$, where $\sigma$ is the unique critical 0 cell, with the connected components that correspond to $y$'s other child nodes. Let $\tau$ be one of these merging 1 cells. Then one of $\tau$'s boundary 0 cells $\sigma_0$ belongs to $X_{M(f)_{y}-\varepsilon}[x]$ and there is a gradient path $\gamma$ in $X_{M(f)_{y}}[\tau]$ that begins with $\tau\supset \sigma_0$. Since maximal gradient paths lead to critical cells and cannot lead through cells of increasing dimension, $\gamma$ has to lead to $\sigma$. If $f$ is weakly Morse--Benedetti, then each inner node $y$ corresponds to a unique merging 1 cell $\tau$ and because, by assumption, only one of $\tau$'s boundary 0 cells $\sigma_0$ belongs to $X_{M(f)_{y}-\varepsilon}[x]$, the gradient path starting at $\tau \supset \sigma_0$ is unique because it is contained in the 1 skeleton of $X$.
        \item  The hyperplanes of $\mathcal{O}(M(f))$ correspond by definition to inclusions of connected components sublevel complexes of $f$. Similar to the proof of (1), cover relations $x\prec y$ of nodes in $M(f)$ correspond to gradient paths, which lead from the merging 1 cells that correspond to $y$ to any critical 0 cell of the connected component that corresponds to $x$. If $x\leq y$ is not a cover relation, then the statement follows inductively because the interval $[x,y]$ is a chain due to $P$ being of merge tree type.  
        \item Since posets of merge tree type have a unique maximum, namely the root, by definition, all elements of $P$ are related by zig-zags of comparison relations. It follows from (2) that cells that correspond to elements of $P$ are related by zig-zags of zig-zags of gradient paths, which are zig-zags of gradient paths. 
    \end{enumerate}
\end{proof}
\begin{Rem}
    Since discrete gradient paths are sequences of matched cells, gradient paths correspond to intersections of hyperplanes of the Morse arrangement $\mathcal{A}(X)$. Even though this observation makes the correspondence between families of gradient paths and hyperplanes in \Cref{PropGradientPahtsAndHyperplanes} (1) into a correspondence between certain intersections of hyperplanes in $\mathcal{A}(X)$ and hyperplanes in $\mathcal{A}(M(f))$, that does not in general lead to a correspondence between cancellations of gradient paths and removals of leaves. 
\end{Rem}
\begin{Cor}\label{Corollary: CancellationsInduceRemovalofLeaves}
    If $f\colon X \rightarrow \mathbb{R}$ is a weak Morse--Benedetti function, then cancellations of critical 0 cells and 1 cells are in 1-1 correspondence with removals of leaves of the induced merge tree. If $f$ is not weakly Morse--Benedetti, cancellations of critical 0 cells and 1 cells still induce removals of leaves of the induced merge tree but the canceled pair is not necessarily unique. 
\end{Cor}
\begin{Def}
    Let $f$ be a discrete Morse function and let $\eta_f$ be a composition of splittings of inner nodes and additions of leaves. Then we denote by $\operatorname{Can}(\eta_f)$ the collection of pairs of simplices whose cancellation via \Cref{Corollary: CancellationsInduceRemovalofLeaves} corresponds to the removal of leaves that $\eta$ conducts. Moreover, we denote by $\operatorname{Mov}(\eta)$ the collection of critical 1-cells whose function values are changed in order to induce the splitting described by $\eta$, as is explained in \Cref{LemmaMoveCrossing}.
\end{Def}

\begin{Def}\label{definitionspaceofmergetrees}
In order to topologize the space of merge trees $Mer$, we introduce the \D{edit path distance}:
Let $\theta\colon P \rightarrow \mathbb{R},\theta'\colon P' \rightarrow \mathbb{R}$ be two merge trees and let $\gamma =(\eta_n\circ\dots\circ \eta_1 \colon s(\theta)  \rightarrow s(\theta'))$ be a sequence of compositions of elementary edit moves in the same direction from $\theta$ to $\theta'$. We define $\gamma_0 \coloneqq \theta$, $\gamma_{n+1}\coloneqq \theta'$, and 
$P_i \coloneqq \operatorname{im}(\eta_i)$. Moreover, we consider sequences of height functions $\bar{\gamma}\coloneqq \gamma_i\colon P_i \rightarrow\mathbb{R}$.
We define the \D{length} of $\gamma$ by
$$
\lambda(\gamma) \coloneqq \inf\limits_{\bar{\gamma}}\sum\limits_{i=0}^{n}\lvert \lvert {\eta_i}_*\gamma_i - \gamma_{i+1}\rvert \rvert ,
$$
where $\lvert \lvert \cdot \rvert \rvert$ denotes the Euclidean norm.
We then define the \D{edit path distance} between $\theta$ and $\theta'$ by 
\begin{align}\label{FormulaEuclideanEditDistance}
\operatorname{d}(\theta,\theta')\coloneqq \inf\limits_{\gamma\colon s(\theta)\rightarrow s(\theta')}\lambda(\gamma) .
\end{align}
 
\end{Def}

\begin{Prop}\label{PropEuclideanEditDistanceIsMetric}
The edit path distance $\operatorname{d}$ from \Cref{definitionspaceofmergetrees} is a pseudo-metric on $Mer$, $Mer_<$, and $Mer_{wb}$. 
\end{Prop}
\begin{proof}
    The equality $\operatorname{d}(\theta,\theta)=0$ holds because the identity is an edit move. Symmetry holds because every edit move has an inverse: adding leaves is inverse to removing them and splitting inner nodes in inverse to merging them. Sequences of edit morphisms can be inverted stepwise by inverting elementary edit morphisms. The triangle inequality follows by construction as a path metric. Moreover, any two merge trees are connected by a sequence of edit moves, for example by merging all inner nodes and removing all leaves, which results in the trivial rooted tree.  
\end{proof}
\begin{Rem}
    The edit path distance cannot be a metric on $Mer$ because of the possibility of constant functions. These can result in paths of length 0 between merge trees that only differ in adjacent nodes with the same function value.\par 
    Moreover, the edit path distance cannot be a metric on $Mer_<$ and $Mer_{wb}$ either due to symmetries of the underlying combinatorial merge trees. It is possible to fix the latter issue by consideration of ordered underlying combinatorial merge trees, i.e.~total order on the sets of children for each inner node. Then it is possible to fix representatives by demanding a suitable notion of compatibility between the orders on the sets of children and function values on the nodes. 
\end{Rem}
\begin{Theo}\label{TheoInducedMergeTreeIsContinuous}
    Let $X$ be a regular CW complex. Then the map $M \colon \mathcal{M}(X)\rightarrow Mer$, where $\mathcal{M}(X)$ is equipped with the path distance\footnote{See \Cref{DefinitionPathDistance}.} and Mer is equipped with the edit path distance, from \Cref{inducedmergetreemap} that maps a discrete Morse function to its induced merge tree is Lipschitz continuous.
\end{Theo}
\begin{proof}
     We have to show that there is a constant $c \in \mathbb{R}$ such that for all $f,g \in \mathcal{M}(X)$ we have $\operatorname{d}(M(f),M(g))\leq c \cdot \operatorname{d}(f,g)$. \par 
     If $f,g$ lie in the interior of the same region $R$ of $\mathcal{H}$, then $f$ and $g$ have the same critical cells, the same gradient field, and induce the same combinatorial merge tree, so $M$ maps $f$ and $g$ to the same subspace, so $M_{\lvert R}$ is just an orthogonal projection, which is Lipschitz continuous with $c=1$. In this case, the path metric and the edit path distance both reduce to the Euclidean metric. \par
     For the case that $f,g$ lie in different regions $f\in R_1, g\in R_2$ of $\mathcal{H}$, we consider that
    by \Cref{PropGradientPahtsAndHyperplanes}, \Cref{Lemmacancelcrossing}, \Cref{LemmaMoveCrossing}, and \Cref{Corollary: CancellationsInduceRemovalofLeaves} the necessary edit morphisms $\tilde{\gamma}$ between the induced combinatorial merge trees correspond to paths $\gamma$ in the space of discrete Morse functions that cross certain hyperplanes $H_i$ of $\mathcal{H}$. \par
    Let $\tilde{\gamma}$ be any sequence of edit moves from $s(M(f))$ to $s(M(g))$ and let $\gamma$ be a corresponding path in the space of discrete Morse functions from $f$ to $g$. We will prove the statement of the theorem by showing that $\lambda(\gamma)\leq(\lambda(\tilde{\gamma}))$.\par
    Let $(\tilde{\gamma}_i)$ be a sequence of intermediate height functions as in \Cref{definitionspaceofmergetrees}. Due to each of the $\tilde{\gamma}_i$ lying inside an intersection of hyperplanes of the braid arrangement on $Mer$ (or in a face reached by exiting a certain intersection of hyperplanes), the corresponding path $\gamma$ has a corresponding subdivision with intermediate functions $h_i$ inside the corresponding intersections of hyperplanes of the braid arrangement (or in a face reached by exiting the corresponding intersection of hyperplanes) on the space of discrete Morse functions. The nodes of the induced merge tree only correspond to certain critical 1- and all critical 0-cells of $X$, which implies for the computation of the length of path segments in the space of merge trees, that the relevant vectors can be embedded as subvectors in the relevant vectors in the space of discrete Morse functions. Hence, we have $\lvert\lvert {\eta_i}_* \gamma_i-\gamma_{i+1}\rvert \rvert \leq \lvert\lvert h_i-h_{i+1}$ for all $i$ even in the space of discrete functions. Moreover, the path $\gamma$ may have to cross additional intersections of hyperplanes in order to stay within the space of discrete Morse functions, which could only potentially increase the length of $\gamma$ due to the triangle inequality. \par
    Since this inequality between the lengths of path segments holds for all segments of each potential path, it also holds for the respective infima and, hence, also for the distance $d(M(f),M(g)\leq d(f,g)$. Hence, the induced merge tree $M$ satisfies Lipschitz-continuity with $c=1$.

\end{proof}
\begin{Rem}
Lipschitz continuity of invariants of filtered spaces is an especially relevant property for applications in TDA. In this framework, being Lipschitz continuous is also called satisfying stability and is a necessary property for invariants to be used for statistical analysis. Stability of merge trees is by no means a new result, for example in \cite{cardona_et_al:LIPIcs.SoCG.2022.24}, the authors have proved stability of the induced merge tree with respect to a large class of distances on merge trees.
\end{Rem}
\begin{Ex}
    To illustrate the ideas of the proof, we give an example of the distance between two discrete Morse functions on a regular CW complex and their induced merge trees:
    \begin{center}
        \begin{tikzpicture}
        \draw (1,-1) node {$f\colon$};
            \draw[fill=lightgray] (5,-1) -- (6.5,-2) -- (6.5,0.25) -- cycle;
    \draw[color=pink] (5,-1) -- (6.5,-2);
    \draw[color=lime] (6.5,-2) -- (6.5,0.25);
	\shade[ball color = gray!40, opacity = 0.4] (3.5,-1) circle (1.5cm);
	\draw (3.5,-1) circle (1.5cm);
	\draw[color=orange] (2,-1) arc (180:360:1.5 and 0.6);
	\draw[dashed,color=orange] (5,-1) arc (0:180:1.5 and 0.6);
    \draw (2,-1) node {$\bullet$};
	\draw (5,-1) node[color=blue] {$\bullet$};
    \node[circle, fill=green, inner sep=1.5pt, outer sep=0pt] (0) at (6.5,-2) {};
    \node[circle, fill=purple, inner sep=1.5pt, outer sep=0pt] (1) at (6.5,0.25) {};
	\draw (2,-1) node[anchor=east] {0};
    \draw (3.5,-1) node[anchor=north] {2};
    \draw (3.5,-1) node[anchor=south] {2};
    \draw (3.5,0.5) node[anchor=south] {3};
    \draw (3.5,-2.5) node[anchor=north] {4};
    \draw (4.9,-1.1) node[anchor=north west] {1};
    \draw (5.75,-0.375) node[anchor=south] {5};
    \draw (5.75,-1.625) node[anchor=north] {4};
    \draw (5.8,-1) node[] {5};
    \draw (6.5,0.25) node[anchor=south] {0};
    \draw (6.5,-2) node[anchor=north] {1};
    \draw (6.5,-1) node[anchor=west] {3};
    \draw (7.5,-1) node {$\xrightarrow[]{M}$};
    \node[circle, fill=pink, inner sep=1.5pt, outer sep=0pt] (a) at (10,0.25) {};
    \draw (10,0.25) node[anchor=south] {4};
    \node[circle, fill=orange, inner sep=1.5pt, outer sep=0pt] (b) at (9,-0.75) {};
    \draw (9,-0.75) node[anchor=east] {2};
    \node[circle, fill=lime, inner sep=1.5pt, outer sep=0pt] (c) at (11,-0.75) {};
    \draw (11,-0.75) node[anchor=west] {3};
    \node[circle, fill=black, inner sep=1.5pt, outer sep=0pt] (d) at (8.5,-1.75) {};
    \draw (8.5,-1.75) node[anchor=north] {0};
    \node[circle, fill=blue, inner sep=1.5pt, outer sep=0pt] (e) at (9.5,-1.75) {};
    \draw (9.5,-1.75) node[anchor=north] {1};
    \node[circle, fill=green, inner sep=1.5pt, outer sep=0pt] (f) at (10.5,-1.75) {};
    \draw (10.5,-1.75) node[anchor=north] {1};
    \node[circle, fill=purple, inner sep=1.5pt, outer sep=0pt] (g) at (11.5,-1.75) {};
    \draw (11.5,-1.75) node[anchor=north] {0};
    \draw (b) -- (a);
    \draw (c) -- (a);
    \draw (b) -- (d);
    \draw (b) -- (e);
    \draw (c) -- (f);
    \draw (c) -- (g);
        \end{tikzpicture}
        \begin{tikzpicture}
        \draw (1,-1) node {$g\colon$};
            \draw[fill=lightgray] (5,-1) -- (6.5,-2) -- (6.5,0.25) -- cycle;
    \draw[color=orange] (5,-1) -- (6.5,-2);
    \draw[color=purple] (6.5,-2) -- (6.5,0.25);
	\shade[ball color = gray!40, opacity = 0.4] (3.5,-1) circle (1.5cm);
	\draw (3.5,-1) circle (1.5cm);
	\draw[color=orange] (2,-1) arc (180:360:1.5 and 0.6);
	\draw[dashed,color=black] (5,-1) arc (0:180:1.5 and 0.6);
    \draw (2,-1) node {$\bullet$};
	\draw (5,-1) node[color=blue] {$\bullet$};
    \node[circle, fill=green, inner sep=1.5pt, outer sep=0pt] (0) at (6.5,-2) {};
    \node[circle, fill=purple, inner sep=1.5pt, outer sep=0pt] (1) at (6.5,0.25) {};
	\draw (2,-1) node[anchor=east] {0};
    \draw (3.5,-1) node[anchor=north] {3};
    \draw (3.5,-1) node[anchor=south] {4};
    \draw (3.5,0.5) node[anchor=south] {4};
    \draw (3.5,-2.5) node[anchor=north] {5};
    \draw (4.9,-1.1) node[anchor=north west] {1};
    \draw (5.75,-0.375) node[anchor=south] {5};
    \draw (5.75,-1.625) node[anchor=north] {3};
    \draw (5.8,-1) node[] {5};
    \draw (6.5,0.25) node[anchor=south] {0};
    \draw (6.5,-2) node[anchor=north] {1};
    \draw (6.5,-1) node[anchor=west] {4};
    \draw (7.5,-1) node {$\xrightarrow[]{M}$};
    \node[circle, fill=purple, inner sep=1.5pt, outer sep=0pt] (a) at (10,0.25) {};
    \draw (10,0.25) node[anchor=south] {4};
    \node[circle, fill=orange, inner sep=1.5pt, outer sep=0pt] (b) at (9,-0.75) {};
    \draw (9,-0.75) node[anchor=east] {3};
    \node[circle, fill=purple, inner sep=1.5pt, outer sep=0pt] (c) at (11,-0.75) {};
    \draw (11,-0.75) node[anchor=west] {0};
    \node[circle, fill=black, inner sep=1.5pt, outer sep=0pt] (d) at (8.5,-1.75) {};
    \draw (8.5,-1.75) node[anchor=north] {0};
    \node[circle, fill=blue, inner sep=1.5pt, outer sep=0pt] (e) at (9.5,-1.75) {};
    \draw (9.5,-1.75) node[anchor=north] {1};
    \node[circle, fill=green, inner sep=1.5pt, outer sep=0pt] (f) at (9,-1.75) {};
    \draw (9,-1.75) node[anchor=north] {1};
    \draw (b) -- (a);
    \draw (c) -- (a);
    \draw (b) -- (d);
    \draw (b) -- (e);
    \draw (b) -- (f);
        \end{tikzpicture}
    \end{center}
    On the side of the induced merge trees, we describe the necessary edit moves from $M(f)$ to $M(g)$: First, we merge the lime-colored inner node labeled 3 with the pink-colored inner node labeled 4. Then we split off a new inner node from the pink node the other way around, label it 3, and merge the orange node with new newly created node. Throughout the whole process, we assign the green node labeled 1 as a child to the inner node we move around. As a result, the fixed nodes have distance zero and the distance between the merge trees is purely given by the edit costs.
    $$d(M(f),M(g))=1+1+1=3.$$
    Regarding the distance between  the discrete Morse functions, we consider that there are 6 cells whose function values with respect to $f$ and $g$ lie in different faces of $\mathcal{H}$: all 1 cells except for the one labeled 5, and the two 2 cells of the sphere on the left. We denote them as follows:
    \begin{center}
        \begin{tikzpicture}
            \draw[fill=lightgray] (5,-1) -- (6.5,-2) -- (6.5,0.25) -- cycle;
    \draw[color=black] (5,-1) -- (6.5,-2);
    \draw[color=black] (6.5,-2) -- (6.5,0.25);
	\shade[ball color = gray!40, opacity = 0.4] (3.5,-1) circle (1.5cm);
	\draw (3.5,-1) circle (1.5cm);
	\draw[color=black] (2,-1) arc (180:360:1.5 and 0.6);
	\draw[dashed,color=black] (5,-1) arc (0:180:1.5 and 0.6);
    \draw (2,-1) node {$\bullet$};
	\draw (5,-1) node[color=black] {$\bullet$};
    \node[circle, fill=black, inner sep=1.5pt, outer sep=0pt] (0) at (6.5,-2) {};
    \node[circle, fill=black, inner sep=1.5pt, outer sep=0pt] (1) at (6.5,0.25) {};
	\draw (2,-1) node[anchor=east] {};
    \draw (3.5,-1) node[anchor=north] {$\gamma$};
    \draw (3.5,-1) node[anchor=south] {$\beta$};
    \draw (3.5,0.5) node[anchor=south] {$\alpha$};
    \draw (3.5,-2.5) node[anchor=north] {$\delta$};
    \draw (4.9,-1.1) node[anchor=north west] {};
    \draw (5.75,-0.375) node[anchor=south] {};
    \draw (5.75,-1.625) node[anchor=north] {$\epsilon$};
    \draw (5.8,-1) node[] {};
    \draw (6.5,0.25) node[anchor=south] {};
    \draw (6.5,-2) node[anchor=north] {};
    \draw (6.5,-1) node[anchor=west] {$\rho$};
        \end{tikzpicture}
    \end{center}
    We consider the relevant differences in function values that create the necessity of crossing certain hyperplanes:
    \begin{multicols}{2}
    \begin{enumerate}
        \item $(\alpha) >f(\beta), \  g(\alpha)=g(\beta)$
        \item $f(\alpha) < f(\epsilon) ,\ g(\alpha) > g(\epsilon)$
        \item $f(\beta)=f(\gamma),\ g(\beta)>g(\gamma) $ 
        \item $f(\beta)<f(\epsilon), \ g(\beta)>g(\epsilon) $ 
        \item $f(\beta)<f(\rho), \ g(\beta)=g(\rho) $
        \item $f(\gamma)<f(\epsilon), \ g(\gamma)=g(\epsilon)$ 
        \item $f(\delta) =f(\epsilon), \ g(\delta) > g(\epsilon)$
        \item $f(\epsilon)>f(\rho), \ g(\epsilon) < g(\rho)$
    \end{enumerate}
\end{multicols}
    It follows that in order to get from $f$ to $g$ in the space of discrete Morse functions, at least 11 entries/exits of hyperplanes need to happen, which result in 11 crossing costs. The crossings related to the changes in (8) and (6) correspond to the above described edit moves on the induced merge trees and have the same cost. The other necessary crossings add to the distance between the functions but are not reflected in the distance between the induced merge tree.\par
    The optimal path from $f$ to $g$ turns out to be given by the direct line segment between $f$ and $g$. If we parametrize it by $\gamma(t)\coloneqq (1-t)\cdot f + t \cdot g$, the hyperplanes that correspond to (3) and (7) are exited at $t=0$, the hyperplanes that correspond to (2) and (8) are being crossed at $t=\frac{1}{2}$, the hyperplane that corresponds to (4) is being crossed at $t=\frac{2}{3}$, and the hyperplanes that correspond to (1), (5), and (6) are entered at $t=1$. The crossing costs are minimized by taking the cut points in the middle between the crossings. This results in a distance of 
    $$
    d(f,g)=\frac{3}{4}+\frac{3}{4}+\frac{1}{4}+\frac{1}{4}+\frac{1}{2}+\frac{1}{2}=3
    $$
    Since in this example all intersections of hyperplanes that are being crossed by the direct line between $f$ and $g$ are actually contained in the space of discrete Morse functions and the function values agree on matched cells, the distance agrees with the one between the induced merge trees. 
\end{Ex}
\subsection{The space of barcodes}
We identify sets of barcodes with sets of order-preserving maps $\operatorname{Map}_\leq(P,\mathbb{R})$ for posets $P$ of barcode type. 
\begin{Def}
    We define the \D{set of barcodes} by 
    $$Bar\coloneqq \coprod\limits_{P \text{ of barcode type}} \operatorname{Map}_\leq(P,\mathbb{R}).$$ 
\end{Def}
\begin{Rem}
    Technically, we have to solve the same set theoretic issues as for the space of merge trees again for the space of barcodes. They can be solved in the same way as for merge trees.
\end{Rem}
\begin{Def}
    Let $P,P'$ be two posets of barcode type. 
    \begin{enumerate}
        \item If $P\cong P'$ as posets, we call an isomorphism $\eta\colon P \rightarrow P'$ a \D{reordering} of the bars. We call $\operatorname{id\colon P \rightarrow P}$ the \D{trivial reordering}.
        \item We say that $P'$ arises from $P$ by \D{collapsing a bar} if $P'$ has one bar fewer than $P$. In that case, a \D{collapse} of a bar is an injective poset map $\eta\colon P'\rightarrow P$. Moreover, we say that $P$ arises from $P'$ by \D{creation} of a new bar. 
    \end{enumerate}
    We refer to all of the mentioned operations as \D{edit moves}. An \D{elementary edit morphism} between posets of barcode type is a composition of an edit move with a reordering of the bars. An \D{edit morphism} is a composition of elementary edit morphisms. \par
    Associated to an edit morphism $\eta\colon P \rightarrow P'$, we define a map $\eta_* \colon \operatorname{Map}_\leq(P,\mathbb{R})\times \operatorname{Map}_\leq(P',\mathbb{R}) \rightarrow \operatorname{Map}_\leq(P',\mathbb{R})$ by 
    $$\eta_*(\beta,\beta')_x\coloneqq \begin{cases}
        \beta_x \text{ if } x \in P,\\
        \beta'_x \text{ if } x \notin P
    \end{cases}.$$
\end{Def} 
\begin{Def}
    The \D{edit path distance} on $Bar$ is defined as follows:\par
    Let $P,P'$ be posets of barcode type. Let $\beta_1 \in \operatorname{Map}_\leq(P,\mathbb{R}), \beta_2 \in \operatorname{Map}_\leq(P',\mathbb{R})$. If $P'$ arises from $P$ by an edit morphism $\eta \colon P\rightarrow P'$, we define $ \operatorname{d}(\beta_1,\beta_2)\coloneqq \operatorname{d^{euc}(\eta_*(\beta,\beta_2),\beta_2)}$
    In the general case, we define the \D{edit path distance} by $$\operatorname{d}(\beta,\beta')\coloneqq \inf\limits_{\sigma}\{\sum\limits_i\operatorname{d}(\beta_i,\beta_{i+1})\},$$where $\sigma =(\beta=\beta_1,\dots, \beta_n=\beta')$ are sequences of barcodes related by edit morphisms.\par    
\end{Def}
\begin{Prop}
    The edit path distance is a pseudo-metric on $Bar$
\end{Prop}
\begin{proof}
    The proof is analogous to the one of \Cref{PropEuclideanEditDistanceIsMetric}.
\end{proof}
As for merge trees, the edit path distance on barcodes cannot be a metric due to symmetries. Again, it is possible to deal with the symmetries by using ordered barcodes to choose specific representatives.\par
We topologize the space of barcodes with the edit path distance.
\begin{Def}
    Let $P$ be a poset of barcode type. We define the \D{order arrangement} $\mathcal{\mathcal{O}(P)}$ in $\operatorname{Map}_\leq(P,\mathbb{R})$ by $\mathcal{O}(P)\coloneqq \{H_x^y \vert x\leq y\}$. \par
\end{Def}
In order to define a map from the space of well-branched merge trees to the space of barcodes, we apply the combinatorial elder rule, see \cite[Definition 2.15]{CURRY2024102031}
\begin{Def}\label{DefInducedBarcode}
    We define a map $B\colon Mer_{wb} \rightarrow Bar$ as follows: \par
    Let $\theta \in Mer_{wb}$ be a well-branched merge tree with an underlying combinatorial merge tree $P$. We define a poset $\tilde{B}(\theta)$ inductively: Add for each inner node $y$ of $P$ an interval $[x,y]$ to $\tilde{B}(\theta)$, where $x$ denotes the leaf of $P_{\leq y}$ with the minimal value $\theta_x$. For all leaf nodes $x$ that are not matched by the rule above, we add an interval $[x,y]$ to $\tilde{B}(\theta)$ where $y$ is $x$'s parent node.\par
    By construction, we have a canonical map $i \colon \tilde{B}(\theta)\rightarrow P$ given by mapping the elements of $\tilde{B}(\theta)$ to their corresponding nodes in $P$. We define $B(\theta)\coloneqq i^*(\theta)=\theta \circ i$. 
\end{Def}
\begin{Prop}\label{PropMapInducedBarcode}
    The map $B \colon Mer_{wb} \rightarrow Bar$ from \Cref{DefInducedBarcode} is Lipschitz continuous with constant $c=1$. 
\end{Prop}
\begin{proof}
Let $\theta$ and $\theta'$ be well-branched merge trees and let $B(\theta),B(\theta')$ be their induced barcodes. Any sequence of edit moves from $\theta$ to $\theta'$ induces a sequence of edit moves from $B(\theta)$ to $B(\theta')$ such that merges/splittings of inner nodes induce reorderings and additions/removals of leaves induced creations/collapses. By definition of the edit costs, additions/removals of leaves have the same cost as their induced creation/collapse of bars. Reorderings can be done simultaneously whereas splittings and merges of inner nodes need to be done consecutively, which is why the distance may decrease when proceeding to the induced barcode.

\end{proof}

\section{Future directions}
The main contribution of this article is a novel framework for the space of discrete Morse functions, which turns out to naturally fit into existing related mathematical concepts. Due to this, there are several potential future directions connected to a deeper investigation of any of the mentioned connections to other mathematical concepts.\par
As briefly mentioned at the end of section 3, a globalization of the concept of the space of discrete Morse functions to a simple homotopy type may be useful for further studies towards the relationship to the smooth case. One could try to approach such studies in a local way, by studying paths in the space of discrete Morse functions on a CW complex $X$ and allowing them to jump to the space of discrete Morse functions on a simple homotopy equivalent CW complex $X'$ when they reach a relevant intersection of hyperplanes, similar to the way we have treated merge trees in this work. However, it is not clear how such an approach may lead to a global understanding of the space of discrete Morse functions on a simple homotopy type. For this, it seems necessary to find a coherent way of gluing spaces of discrete Morse functions on simple homotopy equivalent complexes in a colimit-type construction. Given such a notion of space of discrete Morse functions on a simple homotopy type, one could further study the relationship between smooth and discrete Morse theory by studying an extension of Cerf's map $\eta$ that we extended in \Cref{CorCerf}.\par
In a different direction, one could further study complexes of discrete Morse matchings as sub-posets of the intersection lattice of the Morse arrangement. This may lead to a better understanding of how cancellations of critical cells affect the sublevel filtration induced by the discrete Morse function, which in turn could be used for Morse-theoretic algorithms in topological data analysis.\par
In a similar direction, the relationship to the spaces of merge trees and barcodes provide a novel framework for studying inverse problems in persistence. 
\bibliographystyle{alpha}
\bibliography{References}
\end{document}